\documentclass[11pt]{amsart}
\title{Planar Algebra of the Subgroup-Subfactor}

\author{Ved Prakash Gupta}
\address{The Institute of Mathematical Sciences, Chennai, India}
\email{vpgupta@imsc.res.in}



\usepackage{amsmath}
\usepackage{amsfonts}
\usepackage{latexsym}
\usepackage{amssymb}
\usepackage{amscd}
\usepackage[all]{xy}
\usepackage{epsfig}
\usepackage{graphicx}
\usepackage{hyperref}
\usepackage{eufrak}
\usepackage{psfrag}

\newtheorem{thm}{Theorem}[section]
\newtheorem{lem}[thm]{Lemma}
\newtheorem{cor}[thm]{Corollary}
\newtheorem{prop}[thm]{Proposition}
\newtheorem{defn}[thm]{Definition}
\newtheorem{rem}[thm]{Remark}

\newenvironment{pf}{\noindent{\em 
Proof:}}{\hfill $\Box$}

\begin{document}
\maketitle \hspace*{5mm}{\em Accepted and to appear in Proc.~Indian
  Acad.~Sci.~(Math.~Sci.)}.
\begin{abstract}
  We give an identification between the planar algebra of the
  subgroup-subfactor $R \rtimes H \subset R \rtimes G$ and the
  $G$-invariant planar subalgebra of the planar algebra of the
  bipartite graph $\star_n$, where $n = [G : H]$. The crucial step in
  this identification is an exhibition of a model for the basic
  construction tower, and thereafter of the standard invariant, of $R
  \rtimes H \subset R \rtimes G$ in terms of operator matrices.

  We also obtain an identification between the planar algebra of the
  fixed algebra subfactor $R^G \subset R^H$ and the $G$-invariant
  planar subalgebra of the planar algebra of the `flip' of $\star_n $.

\end{abstract}

\section{Introduction} 
For every pair $H \subset G$ of finite groups, and an outer action
$\alpha$ of $G$ on the hyperfinite $II_1$-factor $R$, we have a dual
pair of subfactors $R \rtimes_{\alpha \slash_H} H \subset R
\rtimes_{\alpha} G$ (the {\em subgroup-subfactor}) and fixed algebra
subfactor $R^G \subset R^H$.

On the other hand, Jones \cite{Jon00} associates a planar algebra
$P(\Gamma)$ to a finite bipartite graph $\Gamma$ with a spin function
$\mu$.

We show that the planar algebras of the subgroup and fixed algebra
subfactors\footnote{Our subfactors are all of type $II_1$.} are
isomorphic to the invariants with respect to an action 
by $G$ as planar algebra automorphisms of $P(\Gamma)$ and
$P(\overline{\Gamma})$, respectively (see Theorem \ref{main-theorem}
and Corollary \ref{dual-theorem}), where $\Gamma$ is the bipartite
graph $\star_n$ ($n = [G : H]$), and $\overline{\Gamma}$ is its `flip'
(with `even' and `odd' vertices interchanged).

We begin $\S2$ with a discussion of matrix algebras, matrix maps and
our notational conventions for them, and go on with certain facts
related to the matrix functor in the context of subfactors.  $\S3$ is
mainly devoted to the exhibition of an explicit model for the basic
construction tower of $R \rtimes H \subset R \rtimes G$ and then its
relative commutants in terms of that model.  In $\S4$, we briefly
recall some basic aspects of planar algebras and discuss the notion of
group action on planar algebras and a related result. We then recall
the notion of planar algebra of a bipartite graph as given by Jones
\cite{Jon00}, and discuss how it and its dual behave under the action
of a finite group on the graph. Finally, in $\S5$, we define the
bipartite graph $\star_n$, where $n = [G : H]$, and analyse the planar
algebra $P (\star_n)^G$.  We then complete the proof of the main
theorem via several lemmas; and end the paper with some corollaries.

There are references in literature giving different models for the
basic construction tower and the standard invariant of the
subgroup-subfactor - see, for instance, \cite{KY92} and \cite{BL}.
Apart from that, Bina Bhattacharyya and Zeph Landau \cite{BL} have
given descriptions of planar algebras of an intermediate subfactor and
a subgroup-subfactor is a particular case of it. However, the
description of the planar algebra of a subgroup-subfactor that we
present here is different from theirs, and is more concrete. One
immediate consequence of this description is that, given any pair of
finite groups $H \subset G$ with index $n$, the planar algebra of the
subgroup-subfactor $R \rtimes H \subset R \rtimes G$ is sandwiched in
between the planar algebras $P(\star_n)^{S_n}$ and $P(\star_n)$ - see
Corollay \ref{sandwich}.

\section{Matrix Functor}

Given an algebra $P$ - all our algebras will be over the field
$\mathbb{C}$ - and a (finite) index set $\Lambda$, we write
$M_{\Lambda}(P)$ for the set of matrices with rows and columns indexed
by $\Lambda$, and entries in $P$; this `matrix algebra' is a *-algebra
(resp., a $II_1$ factor) if $P$ is.

An algebra map $\theta : P \rightarrow Q$ gives rise naturally to
the matrix map $M_{\Lambda}(\theta): M_{\Lambda}(P) \rightarrow
M_{\Lambda}(Q)$ 
defined by
\begin{equation}
  \left[M_{\Lambda}(\theta)(A)\right]_{\lambda_1, \lambda_2} = 
  \theta (A_{\lambda_1, \lambda_2}),\ \lambda_1, \lambda_2 \in \Lambda,\,
  A \in M_{\Lambda}(P). 
\end{equation}
For any two index sets $\Lambda$ and $\Gamma$, and an algebra $P$, we
view $M_{\Lambda}(M_{\Gamma}(P))$ as the set of block matrices, whose
blocks are determined by the index set $\Lambda$ and the matrices in
each block are members of $M_{\Gamma}(P)$. Thus we identify it with
the algebra $M_{\Gamma \times \Lambda}(P)$ via the correspondence
\begin{eqnarray}\label{matrix-matrix-algebras}
  M_{\Lambda}(M_{\Gamma}(P)) \ni A & \longmapsto & 
  \widetilde{A} \in M_{\Gamma \times \Lambda}(P), \nonumber\\
  \widetilde{A}_{(\gamma_1, \lambda_1),(\gamma_2, \lambda_2)} & := &  
  (A_{\lambda_1, \lambda_2})_{\gamma_1, \gamma_2}.
\end{eqnarray}
With this convention, we shall interchangeably use $A$ and
$\widetilde{A}$ during calculations.

Thus, given algebras $P$ and $Q$, finite index sets $\Lambda$ and
$\Gamma$, and an algebra map $\theta : P \rightarrow M_{\Gamma}(Q)$,
the matrix map
\[
M_{\Lambda}(\theta) : M_{\Lambda} (P) \rightarrow M_{\Gamma \times
  \Lambda}(Q)
\]
is given by 
\begin{eqnarray}\label{matrix-map} 
[M_{\Lambda}(\theta)(A)]_{(\gamma_1, \lambda_1 ),
  (\gamma_2, \lambda_2)} & = & \theta (A_{\lambda_1, \lambda_2})_{\gamma_1, \gamma_2} 
 =:  \theta_{\gamma_1, \gamma_2}(A_{\lambda_1, \lambda_2}) ,
\end{eqnarray}
$\forall\ (\gamma_i, \lambda_i) \in \Gamma \times \Lambda,\, i = 1,2,\
\forall\ A \in M_{\Lambda} (P)$.

For later reference, we record a few facts - which are mostly
modifications of certain results proved in \cite{JS97}, and whose
proofs we therefore omit:
\begin{prop}\label{M-functor}
  Suppose $N \subset M \subset^{e_1} M_1$ is the basic construction
  for a pair of $II_1$-factors $N \subset M$ (also called a subfactor)
  with $ \tau^{-1} = [M : N] < \infty$.  Then for any finite index set
  $\Lambda$,
\begin{equation}\label{one}
  M_{\Lambda}(N) \subset M_{\Lambda}(M) \subset M_{\Lambda}(M_1)
\end{equation}
is the basic construction for the subfactor $M_{\Lambda}(N) \subset
M_{\Lambda}(M)$ , where
\begin{enumerate}
\item[(i)] the trace preserving conditional expectation
  $E_{M_{\Lambda}(N)}: M_{\Lambda}(M) \rightarrow M_{\Lambda}(N) $ is
  given by the matrix map $M_{\Lambda}(E_N)$; and 
\item[(ii)] the Jones projection which implements the above conditional
  expectation $E_{M_{\Lambda}(N)}$ is given by the diagonal matrix
  $\tilde{e_1} \in M_{\Lambda}(M_1)$, with diagonal entries given by
\[ 
(\widetilde{e}_1)_{\lambda,\, \lambda} = e_1,\ \forall\,\lambda \in \Lambda.
\]
\end{enumerate}
\end{prop}
\begin{lem}\label{useful-amplification}
  Let $N \subset M$ be a subfactor with integer index $n = [M : N]$.
  Let $\{ \lambda_i : i \in I\}$ be an orthonormal
  basis \footnote{$\{\lambda_i\}$ satisfies $E_n (\lambda_i \lambda_j
    ^*) = \delta^i_j,\, \forall i, j \in I.$} for $M \slash N$, where
  $I : = \{1, \ldots, n\}$. Define $ \theta : M \hookrightarrow
  M_I(N)$ by $\theta_{i, j}(x) = E_N(\lambda_i x \lambda_j^*),\,
  \forall\ x \in M,\, i, j \in I $. Then $\theta$ is a unital
  inclusion and we have an isomorphism of towers
\begin{equation}\label{useful-tower}
  \left( N \subset M \stackrel{\theta}{\hookrightarrow}
    M_I(N) \subset M_I(M) \right) \, \cong\, (N \subset M
  \subset M_1 \subset M_2).
\end{equation} 
In particular, the Jones projections $\widetilde{e}_1 \in M_I (N)$ and
$\widetilde{e}_2 \in M_I(M)$, corresponding to $e_1$ and $e_2$ under
the above tower isomorphism, are given by
\begin{eqnarray}
  (\widetilde{e}_1)_{i, j} &=& E_N (\lambda_i) E_N (\lambda_j^*)\, 
  \mbox{and}\label{e1-tilde} \\
  (\widetilde{e}_2)_{i, j} &=& n^{-1} \lambda_i \lambda_j^* , \label{e2-tilde} 
\end{eqnarray}
respectively, $\forall\, i, j \in I$. And the trace preserving
conditional expectation $E_M : M_I(N) \rightarrow M$ is given by
\begin{equation}\label{conditional-expectation}
E_M ((a_{i, j})) = \sum_{i, j \in I} n^{-1} \lambda_i^* a_{i, j} \lambda_j^*,\
 \forall\ (a_{i, j}) \in  M_I(N).
\end{equation}
\end{lem}

\begin{lem}\label{theta-k}
  With $N \subset M$, $\{\lambda_i : i \in I \}$ and $\theta$ as in Lemma
  \ref{useful-amplification}, for each $k \geq 1$, define $\theta
  ^{(k)} : M \rightarrow M_{I^k} (N)$ by 
\begin{eqnarray*}\label{theta-k-equation}
  \theta^{(k)}_{\underline{i},\, \underline{j}}(x) & = & \theta_{i_1,\, j_1}
  (\theta_{i_2,\, j_2}  (\cdots \theta_{i_k,\, j_k}(x) \cdots )) \\ 
  & = & E_N (\lambda_{i_1}E_N( \lambda_{i_2}E_N(\cdots 
  E_N (\lambda_{i_k} x \lambda_{j_k}^* ) 
  \cdots )\lambda_{j_2}^* )\lambda_{j_1}^*),
\end{eqnarray*}
for all $x \in M,\ \underline{i} = (i_1, \ldots, i_k),\,
\underline{j} = (j_1, \ldots, j_k) \in I^k$.  Then $\theta ^{(k)}$ is
a unital $*$-homomorphism for each $k \geq 1$.
\end{lem}

\section{The Subgroup-Subfactor}
In this section, we set up the notation and the model for the
subgroup-subfactor which we shall find convenient to use.

Thus, let $R$ be the hyperfinite $II_1$-factor viewed as sitting in
$\mathcal{L}(L^2 (R))$ as left multiplication operators.  Let $G$ be a
finite group and $\alpha : G \rightarrow Aut(R)$ be an outer action of
$G$ on $R$, i.e., $\alpha_g$ is not inner for all $e \neq g \in G$.

We think of $R \rtimes_{\alpha} G$ as the $II_1$-factor $(R \cup \{u_g
: g \in G\})'' \subset \mathcal{L}(L^2 (R))$, where $u_g (\hat{x}) :=
\widehat{\alpha_g(x)},\, x \in R,\, g \in G$ - see \cite[\S
A.4]{JS97}.

Let $H$ be a subgroup of the group $G$. Once and for all, we fix a set
of representatives $\{g_1, \ldots, g_n\} $ for the right $H$-coset
decomposition of $G$ with $g_1 = e$, i.e., $G = \sqcup_{i = 1}^n H
g_i$, where $n = [G : H]$. Then $\{u_{g_i} : i \in I\}$ is a(n
orthonormal) basis for $R \rtimes G \slash R \rtimes H$, where $I :=
\{1, \ldots, n\}$.

The following result is probably folklore.
\begin{prop}\label{fixed-algebra-subfactor}
\begin{enumerate}
\item $R^G \subset R^H$ is an irreducible subfactor with index $[G :
  H]$ and has finite depth.
\item $(R^G \subset R^H) \, \cong \, ( M \subset M_1 )$, where $M_1$
  is the $II_1$-factor obtained by the basic construction of the
  subgroup-subfactor $N := R \rtimes H \subset R \rtimes G =: M$.
\end{enumerate}
\end{prop}

Taking $N = R \rtimes H \subset R \rtimes G = M$, by Lemma
\ref{useful-amplification}, with respect to the orthonormal basis
$\{u_{g_i} : i \in I \}$ for $M \slash N$, we have
\[
( N \subset M \subset M_1 \subset M_2)\, \cong\, (N \subset M
\stackrel{\theta}{\hookrightarrow} M_I (N) \subset M_I (M) ).
\]
Then, by repeated applications of Proposition \ref{M-functor}, we see
that

\begin{equation}\label{each-stage}
\begin{array}{rcl}
  &  M_{2k-1} \subset M_{2k} \subset M_{2k+1} \subset M_{2k+2}   & \\
  & \cong & \\
  &   M_{I^{k}}(N) \subset M_{I^{k}}(M) \stackrel{\Theta_{k+1}}
  {\hookrightarrow} M_{I^{k+1}}(N) \subset M_{I^{k+1}}(M), &
  \forall\, k \geq 0,
\end{array}
\end{equation}
where, as in $\S 2$, for each $k \geq 1$, we have inductively
identified $M_{I^k}(X)$ with $M_{I} ( M_{I^{k-1}} (X)) $ for $X \in
\{N, M\}$, and $ \Theta_{k+1}:= M_I (\Theta_k)$ with $\Theta_1 :=
\theta$. Thus, we find the following:

\begin{thm}\label{tower-model} With notations as in the preceding
  paragraph, the tower
  \[
  N \subset M \stackrel{\Theta_1}{\hookrightarrow} M_I(N) \subset
  M_I(M) \stackrel{\Theta_2}{\hookrightarrow} \cdots \subset
  M_{I^{k}}(M) \stackrel{\Theta_{k+1}}{\hookrightarrow}
  M_{I^{k+1}}(N)\subset \cdots
  \]
  is a model for the basic construction tower of the
  subgroup-subfactor $R \rtimes H \subset R \rtimes G$.
\end{thm}

Given a $k$-tuple $\, \underline{i} \in I^k$, we write
$\underline{i}_{[r,\, s]}$ for the $(s-r+1)$-tuple $(i_r, i_{r+1},
\ldots, i_s)$ for all $1 \leq r < s \leq k$. In fact, if $r = 1$
(resp., $s = k$), we simply write $\underline{i}_{\ s]}$ (resp.,
$\underline{i}_{[r }$) for $\underline{i}_{[r,\,s]}$. And, for any $j
\in I$, we write $(j,\underline{i}_{[r,\, s]} )$ (resp.,
$(\underline{i}_{[r,\, s]}, j )$) for the tuple $(j, i_r, \ldots,
i_s)$ (resp., $(i_r, \ldots, i_s, j)$).

\begin{cor}\label{model-facts}
  For each $k \geq 1$, let the Jones projections $e_{2k-1} $ and
  $e_{2k}$ be mapped to the operator matrices $\widetilde{e}_{2k-1}
  \in M_{I^k} (N) $ and $\widetilde{e}_{2k} \in M_{I^{k}} (M)$,
  respectively, under the identifications (\ref{each-stage}). Then
\begin{eqnarray}
  (\widetilde{e}_{2k-1})_{\underline{i},\, \underline{j}} & = & 
  \delta^{\underline{i}}_{\underline{j}} 
  \ \delta^{i_1}_{1}, \ 
  \forall\ \underline{i},\, \underline{j} \in I^{k};\ \mbox{and} \\
  (\widetilde{e}_{2k})_{\underline{i},\, \underline{j}} & = & n^{-1}\,
  \delta^{\underline{i}_{[2} }_{\underline{j}_{[2}}\  g_{i_1} g_{j_1}^{-1},\ 
  \forall\   \underline{i},\, \underline{j} \in I^{k}.  
\end{eqnarray}
And, for each $k \geq 0$, the trace preserving conditional
expectations $E_{M_{2k}}: M_{I^{k+1}}(N) \rightarrow M_{I^k}(M)$ and
$E_{M_{2k-1}}: M_{I^{k}}(M) \rightarrow M_{I^k}(N) $ are given by the
matrix maps
\begin{equation}
E_{M_{2k}} = M_{I^k} (E_M)\ \mbox{and}\ E_{M_{2k-1}} = M_{I^k} (E_N),
\end{equation}
respectively, where $E_M$ is given as in Lemma
\ref{useful-amplification}, and as usual\\ $M_{I^{k+1}}(T) :=
M_I(M_{I^{k}}(T))$ for $ T \in \{ E_M, E_N\}$.
\end{cor}
The proof of Corollary \ref{model-facts} relies on Proposition
\ref{M-functor} and Lemma \ref{useful-amplification}, while easy
induction arguments yield the following results.

\begin{lem}\label{Theta_k}
For each $k \geq 1$, we have
\begin{equation}
  \Theta_k ((a_{\underline{i},\, \underline{j}} ))_{\underline{u},\, 
    \underline{v}} = \theta_{u_1, v_1} 
  (a_{\underline{u}_{[2},\, \underline{v}_{[2}}),\, \forall\,
  (a_{\underline{i},\, \underline{j}}) \in M_{I^k} (M),\, 
  \underline{u},\, \underline{v} \in I^{k+1}.
\end{equation}
\end{lem}

\begin{prop} For each $k \geq 1$, we have
\begin{equation}
\Theta_k \circ \cdots \circ \Theta_1 \, = \, \theta^{(k)}.
\end{equation}
\end{prop}
We now recall some notations and facts from \cite{JS97}.

Given a $k$-tuple $\underline{i}= (i_1, \ldots, i_k) \in I^k,\, k \geq
1$, we write $\sqcap g_{\underline{i}}$ for the product $g_{i_1}
g_{i_2} \cdots g_{i_k}$. For each $k \geq 1$, we have an action
$\beta^k$ of $G$ on the set $I^k$: \\
For $\underline{i}, \underline{j} \in I^k$,
\[
  \beta^k_g (\underline{j}) = \underline{i} \Longleftrightarrow
  H g_{j_l}g_{j_{l+1}} 
  \cdots g_{j_k} g^{-1} = H g_{i_l}g_{i_{l+1}} 
  \cdots g_{i_k},\, \forall\, 1 \leq l \leq k. 
\]

For each $k \geq 1$, the relative commutants of $R$ are given by

\begin{eqnarray}
\qquad  \lefteqn{ R' \cap M_{2k-1} \cong \theta^{(k)} (R)' \cap M_{I^k}(N) 
    = \{ (x_{\underline{i},\, \underline{j}}) \in M_{I^k}(N) : } 
  \label{rel-com-1}\\
  &\qquad \exists\ (C_{\underline{i},\, \underline{j}})\in M_{I^k}
  (\mathbb{C}) \ \mbox{such that}\
  x_{\underline{i},\, \underline{j}} = C_{\underline{i},\,
    \underline{j}}\ u_{(\sqcap g_{\underline{i}})(\sqcap
    g_{\underline{j}})^{-1}}& \nonumber\\ 
  &\qquad \quad \mbox{and}\
  C_{\underline{i},\, \underline{j}} = 0\,
  \mbox{unless}\ H\sqcap g_{\underline{i}} = H \sqcap
  g_{\underline{j}},\ \forall\ \underline{i},\, \underline{j} \in I^k
  \};\ \mbox{and}& \nonumber\\
 \qquad \lefteqn{R' \cap M_{2k} \cong \theta^{(k)} (R)' \cap M_{I^k}(M)   =
    \{(x_{\underline{i},\, \underline{j}}) \in M_{I^k}(M) : } 
\label{rel-com-2}\\ 
&  \ \exists\ (C_{\underline{i},\, \underline{j}})\in M_{I^k}(\mathbb{C}) \ \ni 
x_{\underline{i},\, \underline{j}} = C_{\underline{i},\, 
  \underline{j}}\ u_{(\sqcap g_{\underline{i}})(\sqcap g_{\underline{j}})^{-1}}, 
\ \underline{i},\, \underline{j} \in I^k \}.& \nonumber
\end{eqnarray} 

From these descriptions, it immediately follows that the relative
commutants in the grid of the standard invariant of $N \subset M$ are
given by
\begin{eqnarray}
  \lefteqn{N' \cap M_{2k-1} \cong \theta^{(k)} (N)' \cap M_{I^k}(N) 
    = \{ (x_{\underline{i},\, \underline{j}}) 
    \in \theta^{(k)} (R)' \cap M_{I^k}(N) :} \label{rel-com-3} \\
  &\qquad  C_{\underline{i},\, \underline{j}} = 
  C_{\beta^k_h(\underline{i}),\, \beta^k_h(\underline{j})},
  \forall\, h \in H,\,  
  \underline{i},\, \underline{j} \in I^k \};& \nonumber \\ 
  \lefteqn{N' \cap M_{2k} \cong \theta^{(k)} (N)' \cap M_{I^k}(M) 
    = \{ (x_{\underline{i},\, \underline{j}}) 
    \in \theta^{(k)} (R)' \cap M_{I^k}(M) :}\label{rel-com-4}\\
  & \qquad C_{\underline{i},\, \underline{j}} = 
  C_{\beta^k_h(\underline{i}),\, \beta^k_h(\underline{j})},
  \forall\, h \in H,\, 
  \underline{i},\, \underline{j} \in I^k \}; & \nonumber\\
  \lefteqn{M' \cap M_{2k-1} \cong \theta^{(k)} (M)' \cap M_{I^k}(N)  
    = \{(x_{\underline{i},\, \underline{j}}) 
    \in \theta^{(k)} (R)' \cap M_{I^k}(N) :}  \label{rel-com-5}\\
  &\quad \qquad C_{\underline{i},\, \underline{j}} = 
  C_{\beta^k_g(\underline{i}),\, \beta^k_g(\underline{j})},
  \forall\, g \in G,\,   \underline{i},\, \underline{j} \in I^k \};\
  \mbox{and} & \nonumber \\ 
  \lefteqn{M' \cap M_{2k} \cong \theta^{(k)} (M)' \cap M_{I^k}(M) 
    = \{ (x_{\underline{i},\, \underline{j}}) 
    \in \theta^{(k)} (R)' \cap M_{I^k}(M) : }  \label{rel-com-6}\\
  & \qquad  C_{\underline{i},\, \underline{j}} = 
  C_{\beta^k_g(\underline{i}),\, \beta^k_g(\underline{j})},
  \forall\, g \in G,\,   \underline{i},\, \underline{j} \in I^k \}, 
  & \nonumber
\end{eqnarray}
for all $k \geq 1 $, where in each of the above equations the scalar
matrix $(C_{\underline{i},\, \underline{j}})$ corresponds to the
operator matrix $ (x_{\underline{i},\, \underline{j}})$ as in
(\ref{rel-com-1}) and (\ref{rel-com-2}).

For each $k \geq 1$, and given $k$-tuples $\underline{i},\,
\underline{j} \in I^k $, we define operator matrices
$[\underline{i},\, \underline{j}]^{ev} \in M_{I^k}(M) $ and
$[\underline{i},\, \underline{j}]^{od} \in M_{I^k}(N)$ by
\begin{eqnarray}
  [\underline{i},\,\underline{j}]^{ev}_{\underline{u},\,
    \underline{v}} & = & \delta^{\underline{i}}_{\underline{u}} \
  \delta^{\underline{j}}_{\underline{v}}\ u_{(\sqcap
   g_{\underline{i}})(\sqcap g_{\underline{j}})^{-1}}\  \mbox{and} \\
{[\underline{i},\,
 \underline{j}]}^{od}_{\underline{u},\, \underline{v}} 
  & = & \delta^{\underline{i}}_{\underline{u}} \
  \delta^{\underline{j}}_{\underline{v}}\ \ 1_H((\sqcap
  g_{\underline{i}}) (\sqcap g_{\underline{j}})^{-1} )\, u_{(\sqcap
    g_{\underline{i}}) (\sqcap g_{\underline{j}})^{-1}},
\end{eqnarray}
respectively, $ \forall\ \underline{u},\, \underline{v} \in I^k$,
where $1_H (x) := 1$ if $x \in H$ and $0$ otherwise. The superscripts
$ev$ and $od$ refer to the ambient spaces $M_{2k}$ and $M_{2k-1}$,
respectively.

Further, for each $k \geq 1$, we set $Y_k = \{(\underline{i},\,
\underline{j}) \in I^k \times I^k : H \sqcap g_{\underline{i}} = H
\sqcap g_{\underline{j}}\}$. With these simplified notations, it is
not hard to show that, for each $k \geq 1$,
\begin{enumerate}
\item $G$ acts diagonally by $\beta^k $ on $I^k \times I^k$ and $Y_k$
  is invariant under this action, i.e., $ (\beta^k_g
  (\underline{i}),\, \beta^k_g(\underline{j})) \in Y_k,\ \forall\ g
  \in G,\, (\underline{i},\, \underline{j}) \in Y_k $;
\item $\{[\underline{i},\, \underline{j}]^{od} : (\underline{i},\,
  \underline{j}) \in Y_k\}$ forms a basis for $R'\cap M_{2k-1}$; and
\item $\{[\underline{i},\, \underline{j}]^{ev} : (\underline{i},\,
  \underline{j}) \in I^k \times I^k\}$ forms a basis for $R'\cap M_{2k}
  $.
\end{enumerate}
Thus, it makes sense to write $g [\underline{i},\, \underline{j}]^x $
for $ [\beta^k_g(\underline{i}),\, \beta^k_g(\underline{j})]^x$,
$\forall\ g \in G,\ (\underline{i},\, \underline{j}) \in I^k \times
I^k$, $k \geq 1,\ x \in \{ ev,\ od \}$.

Given a $G$-action $\upsilon$ on a set $X$, we write $G\backslash X
(\upsilon)$ for a set of representatives of $G$-orbits of $X$ under
the action $\upsilon$, and when the group action is clear from the
context we simply denote it by $G\backslash X$. For instance, we
write $G \backslash Y_k$ for a set of orbit representatives for the
diagonal $G$-action on $Y_k$ by $\beta^k$.

The descriptions (\ref{rel-com-3}$-$\ref{rel-com-6}) give the
following list of bases for the relative commutants:
\begin{lem}\label{rel-com-bases}
  For each $k \geq 1$, for every choice of orbit representatives $ H
  \backslash X $ and $G\backslash X$ ($X = Y_k,\, I^k \times I^k$),
\begin{enumerate}
\item $\{\sum_{g \in G} g [\underline{i},\, \underline{j}]^{od}:(
  \underline{i},\, \underline{j}) \in  G\backslash Y_k\}$ (resp.,
  $\{\sum_{h \in H} h [\underline{i},\, \underline{j}]^{od}:(
  \underline{i},\, \underline{j})\, \in H \backslash Y_k\}$ ) forms a
  basis for $M' \cap M_{2k-1}$ (resp., $N'\cap M_{2k-1}$).
\item $\{\sum_{g \in G} g [\underline{i},\, \underline{j}]^{ev}: (
  \underline{i},\, \underline{j}) \in  G\backslash (I^k \times
  I^k)\}$ (resp.,$\{\sum_{h \in H} h [\underline{i},\,
  \underline{j}]^{ev}:( \underline{i},\, \underline{j}) \in   H
  \backslash (I^k \times I^k)\}$ ) forms a basis for $M' \cap M_{2k}$
  (resp., $N'\cap M_{2k}$).
\end{enumerate}
\end{lem}
Since the above bases do not vary with the choice of orbit
representatives, we shall no more trouble ourselves by repeating the
phrase ``given a set of representatives of $K$-orbits of $X$'' and
simply write $K \backslash X$ for such a choice.

\section{Planar Algebras}
\subsection{Introduction}
We refer to \cite{KLS03, KS04} for the basic aspects of planar
algebras and thus follow the notations and terminologies therein. The
goal of this section is to analyse the planar algebra associated to a
bipartite graph and its dual, and to see how it responds to an
appropriate action of a finite group on the bipartite graph.

Basically, a planar algebra $P$ is a collection of vector spaces $\{
P_k : k \in Col\}$ admitting an `action' of the operad of coloured
tangles, where $Col : = \{0_{\pm}, 1, 2, \ldots \}$.

\begin{figure}[h]
\begin{center}
\psfrag{0+}{$1^{0_+}$}
\psfrag{0-}{$1^{0_-}$}
\psfrag{k}{$1^k$ ($k$-strings)}
\includegraphics[scale=0.45]{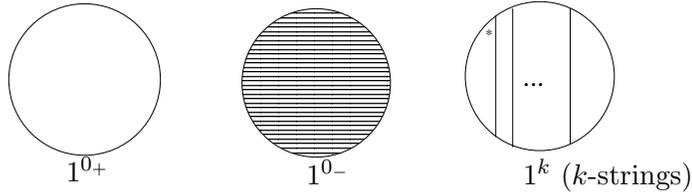}
\end{center}
\caption{The Unit Tangles}\label{unit-tangles}
\end{figure}

\begin{figure}
\begin{center}
\psfrag{1}{$_{^{D_1}}$}
\psfrag{0+}{$I_{0_+}^1$}
\psfrag{0-}{$I_{0_-}^1$}
\psfrag{k}{$I_k^{k+1}$}
\includegraphics[scale=0.40]{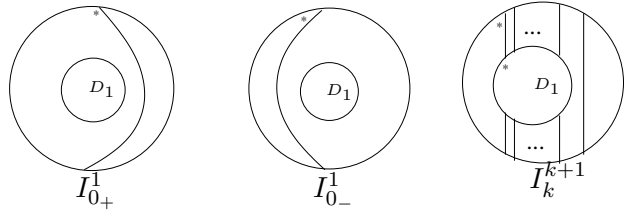}
\end{center}
\caption{The Inclusion Tangles}\label{inclusion-tangles}
\end{figure}

\begin{figure}
\begin{center}
  \psfrag{1}{$_{^{D_1}}$} \psfrag{2}{$_{^{D_2}}$}
  \psfrag{0+}{$M_{0_+}$} \psfrag{0-}{$M_{0_-}$}
  \psfrag{k}{$\begin{array}{c} M_{k} : \ \mbox{all
        three}\\\mbox{are k-discs }\end{array} $}
\includegraphics[scale=0.55]{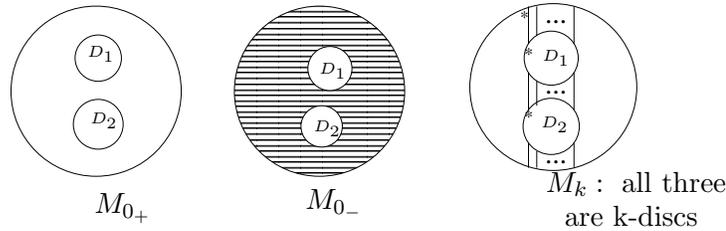}
\end{center}
\caption{The Multiplication Tangles}\label{multiplication-tangles}
\end{figure}

\begin{figure}
\begin{center}
\psfrag{1}{$_{^{D_1}}$}
\psfrag{0+}{$E_{1}^{0_+}$}
\psfrag{0-}{$E_{1}^{0_-}$}
\psfrag{k+1}{$E_{k+1}^{k}$}
\psfrag{1'}{$(E')^1_1$}
\psfrag{k}{$(E')^k_k$}
\includegraphics[scale=0.45]{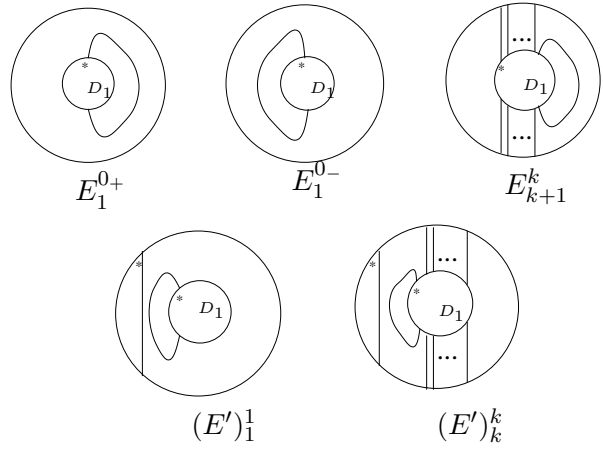}
\end{center}
\caption{The Conditional Expectation
  Tangles}\label{conditional-tangles}
\end{figure}

\begin{figure}
\begin{center}
\psfrag{1}{$_{^{D_1}}$}
\psfrag{R4}{$R_4$}
\psfrag{E1}{$\mathcal{E}^{1+1}$}
\psfrag{Ek}{$\mathcal{E}^{k+1}$}
\includegraphics[scale=0.40]{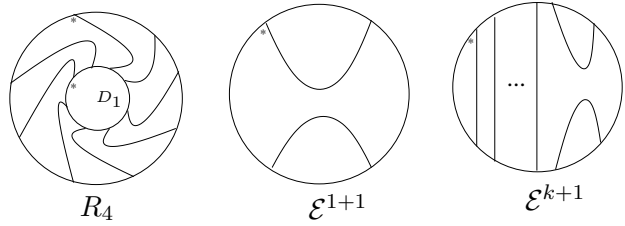}
\end{center}
\caption{The Rotation and  Jones' Projections Tangles }\label{rotn-jones-tangles}
\end{figure}

For the sake of completeness, we list some important coloured tangles
in Figures \ref{unit-tangles}$-$\ref{rotn-jones-tangles}.  With
respect to composition of tangles, \cite{KS04} contains the following
generating sets:
\begin{thm}\label{gen-tangles}\cite{KS04}
Let $\mathcal{T}$ be a collection of coloured tangles containing 
\begin{itemize}\setlength\itemsep{-1mm}
\item $\mathcal{G}_0 := \{ 1^{0_{\pm}}\} \cup \{ E^k_{k+1}, M_k,
  I^{k+1}_{k}:\, k \in\, Col \} \cup \{\mathcal{E}^{k+1}, (E')^k_k : k
  \geq 1\},$ or
\item $\mathcal{G}_1 := \{ 1^{0_{\pm}}\} \cup \{ E^k_{k+1}, M_k,
  I^{k+1}_{k}:\, k \in\, Col \} \cup \{ R_k : k \geq 2\},$
\end{itemize}
and suppose $\mathcal{T}$ is closed under composition of coloured
tangles, whenever it makes sense. Then $\mathcal{T}$ contains all
coloured tangles.
\end{thm}


\begin{defn}\label{planar-morphism} 
  Let $P$ and $Q$ be two planar algebras.  A planar algebra morphism
  from $P$ to $Q$ is a collection $\varphi = \{ \varphi_k : k \in
  Col\}$ of linear maps $\varphi_k : P_k \rightarrow Q_k$ which
  commutes with all the tangle maps, i.e., if $T$ is a $k_0$-tangle
  with $b$ internal boxes of colours $k_1, \ldots, k_b$, respectively,
  then

\begin{equation}
  \varphi_{k_0} \circ Z^P_T 
  = \left\{ 
    \begin{array}{ll}  Z^Q_T \circ (\otimes_{i = 1}^b\,  \varphi_{k_i}) &,\,
      \mbox{if}\ b > 0;\, \mbox{and}  \\
      Z^Q_T&,\, \mbox{if}\ b = 0.\end{array}
  \right.
\end{equation}
$\varphi$ is said to be a planar algebra isomorphism if the maps
$\varphi_k$ are all linear isomorphisms.
\end{defn}
As in the above definition, we shall follow the convention that the
empty tensor product of vector spaces is taken to be the scalars
$\mathbb{C}$.

Repeated applications of the compatibility condition for tangle maps
with respect to composition of tangles give the following useful
result.

\begin{lem}\label{closedness-of-tangles}
  Let $P$ and $Q$ be planar algebras, and $\varphi_k : P_k \rightarrow
  Q_k,\, k \in\, Col$ be linear maps. If $\mathcal{T}$ is the
  collection of those tangles $T$ for which the equation in Definition
  \ref{planar-morphism} holds, then $\mathcal{T}$ is closed under
  composition of tangles.
\end{lem}

We shall primarily deal with subfactor planar algebras.  Basically,
these planar algebras are connected (i.e., $dim P_{0_{\pm}} = 1$),
have positive modulus (i.e., loops come out as constants), are
spherical (i.e., tangle maps for $0_{\pm}$-tangles are isotopy
invariant on the $2$-sphere), their constituent vector spaces are
finite dimensional $C^*$-algebras, and the tangle maps satisfy the
$*$-compatibility condition
\begin{equation}\label{*-condition}
Z_T(x_1 \otimes \cdots \otimes x_b)^* = Z_{T^*}(x_1^* \otimes \cdots
\otimes x_b^*),
\end{equation}
where $T$ is a coloured tangle as in Definition \ref{planar-morphism}
and $ x_i \in P_{k_i}, \, 1 \leq i \leq b$, and $T^*$ is the adjoint
of $T$ - see \cite{KS04}.  Along with these, if the modulus is
$\delta$, then for each $k \geq 0$, the pictorial trace $tr_{k+1} :
P_{k+1} \rightarrow \mathbb{C}$ defined by
\begin{equation}\label{pictorial-trace}
  tr_{k+1} (x) 1_{0_+} = \delta^{-k-1} 
Z^P_{E^{0_+}_{k+1}} (x),\ \forall\   x \in P_{k+1},
\end{equation} 
where $E^{0_+}_{k+1} := E^{0_+}_{1} \circ E^{1}_{2} \circ \cdots \circ
E^{k}_{k+1}$, is a faithful positive trace on $P_{k+1}$.

\begin{lem}\cite{KS06}\label{irreducible-spherical}
  If $P$ is a connected and irreducible\footnote{A planar algebra $P$
    is said to be irreducible if $dim P_1 = 1$.} planar algebra with
  modulus $\delta > 0$, then $P$ is spherical.
\end{lem}

\begin{rem}\label{pictorial-trace-conditional} 
  For a connected and spherical planar algebra $P$ with modulus
  $\delta > 0$, and each $P_k$ being a finite-dimensional
  $C^*$-algebra,
\begin{enumerate}
\item the pictorial traces $tr_m : m \in Col$ are consistent with
  respect to inclusions and thus define a global  trace on
  $P$, where for $0_{\pm}$ the trace for $P_{0_{\pm}} \cong
  \mathbb{C}$ is the obvious identity map; and
\item if $tr_{m}$ is faithful for all $m \geq 1$, then for each $k \in
  Col$, $ \frac{1}{\delta}\, Z^P_{E_{k+1}^k}$ (resp., $
  \frac{1}{\delta}\, Z^P_{(E')^{k+1}_{k+1}}$ ) is the unique
  $tr_{k+1}$ preserving conditional expectation from $P_{k+1}$ onto
  $P_{k}$ (resp., $P_{1,\, k+1}$ ), where $P_{1,\, k+1} := Image
  (Z^P_{(E')_{k+1}^{k+1}})$.
\end{enumerate}
\end{rem}
The importance of planar algebras in subfactor theory lies in the
following theorem of Jones:
\begin{thm}\label{jones}\cite{Jon} For every extremal
  $II_1$-subfactor $N \subset M$ of finite index $\delta^2$, taking
  $P_{0_{\pm}} = \mathbb{C}$ and $P_k = N' \cap M_{k-1},\, k \geq 1$,
  there exists a unique subfactor planar algebra structure on the
  collection $P^{N \subset M} := \{ P_k : k \in \, Col\} =: P$
  satisfying:
\begin{enumerate}
\item $Z^P_{\mathcal{E}^{k+1}} (1) = \delta\, e_k,\, \forall\, k \geq 1;$
\item $Z^P_{(E')^k_k} (x) = \delta\, E_{M' \cap M_{k-1}} (x),\, \forall\,
  x \in N' \cap M_{k-1},\, \forall\, k \geq 1;$
\item $Z^P_{E_{k+1}^k} (x) = \delta\, E_{N'\cap M_{k-1}},\, \forall\, x
  \in N' \cap M_k, \forall\, k \in\, Col$, where for $k=0_{\pm}$, the
  equation is read as 
  \[Z^P_{E_{1}^{0_{\pm}}} (x) = \delta\, tr_{M}(x),\, \forall\, x \in N'
    \cap M.\]
\end{enumerate}
\end{thm}

We now discuss group actions on planar algebras.
\begin{defn}\label{gp-action-plnr-algebra}
  Let $G$ be a finite group and $P = \{ P_k: k \in Col\}$ be a planar
  algebra.  We say that $G$ acts on $P$ if for each $k \in\, Col$, we
  have group homomorphisms $\alpha_k : G \rightarrow GL (P_k)$ such
  that, for each $g \in G$, $\alpha (g) := \{\alpha_k (g): k \in
  Col\}$ is a planar algebra automorphism of $P$.
\end{defn} 
For convenience, we write $g x $ for the element $\alpha_k (g)(x)$,
for $g \in G, x \in P_k$ and $k \in Col$.  Under such an action,
taking 
\[
P^G_k = \{x \in P_k : g x = x,\, \forall g \in G\},\, k \in Col,
\]
we note that the collection $P^G := \{ P^G_k : k \in Col \}$ is a
planar sub-algebra of $P$. Furthermore, the above action induces a
$G$-action on the dual planar algebra $^{-}P$ - see \cite{KS04} - as
well.  Indeed, for each $k > 0$, since $^{-}P_k = P_k$, $G$ acts on
$^{-}P_k $ as it does on $P_k$; and as $^-{P}_{0_{\pm} } = P_{0_{\mp}}
$, again $G$ acts on these as it did on $P_{0_{\mp}}$. The constituent
vector spaces of $^-(P^G)$ and $(^-P)^G$ being the same, the following
is a tautology:

\begin{prop}\label{G-invariant}
  If a finite group $G$ acts on a planar algebra $P$, then with the
  induced $G$-action on $^-P$, the planar algebras $^-(P^G)$ and
  $(^-P)^G$ are isomorphic.
\end{prop}
\subsection{Planar Algebra of the Bipartite
  Graph}\label{pln-bipartite-gph}

Jones \cite{Jon00} associated a planar algebra to a given finite
(possibly with multiple edges), connected bipartite graph with a
``spin function''.  However, in \cite{Jon00}, there is a slight
conflict between the notion of ``state'' and its ``compatibility''
with ``loops''. Thus the planar algebra that we associate here to a
data as above is an appropriate modification of the one given there;
and the results obtained there hold verbatim.

Let $\Gamma = (\mathcal{U}^+,\, \mathcal{U}^-,\, \mathcal{E})$
be a finite, connected bipartite graph with even and odd vertices
$\mathcal{U}^+$ and $\mathcal{U}^-$, respectively, edge set
$\mathcal{E}$ and a {\em spin function} $\mu :\mathcal{U}:=
\mathcal{U}^+ \sqcup \mathcal{U}^- \rightarrow (0, \infty)$.

For each $k >0$, we denote a loop on $\Gamma$ of length $2k$, based at
a vertex $\pi_0 \in \mathcal{U}^+$ by a pair $(\pi, \epsilon)$ of maps
$\pi : \{0, 1, \ldots, 2k-1 \} \rightarrow \mathcal{U}$ and $ \epsilon
: \{ 0, 1, \ldots, 2k-1\} \rightarrow \mathcal{E}$ such that $\pi_{2i}
\in \mathcal{U}^+$ and $\pi_{2i+1} \in \mathcal{U}^-$ for all $0 \leq
i \leq k-1$; and for $0 \leq i \leq 2k-1$, $\epsilon (i)$ is the edge
joining the vertices $\pi_i$ and $\pi_{i+1}$, where we follow, here
and elsewhere, the convention of labelling the vertices of $2k$-loops
{\em modulo $2k$}, thus $\pi_{2k} = \pi_0$.  The following pictorial
form of the loop $(\pi, \epsilon)$ is quite useful:
\[
\xymatrix @-1.1pc{
  & \pi_1 \ar[r]^{\epsilon_1} & \pi_2 \ar[r] & \cdots
  \ar[r] & \pi_{k-1} \ar[rd]^{\epsilon_{k-1}} &  \\
  \pi_0 \ar[ru]^{\epsilon_0} & & & & & \pi_{k} \ar[ld]^{\epsilon_{k}} .\\
  & \pi_{2k-1} \ar[lu]^{\epsilon_{2k-1}} & \pi_{2k-2}
  \ar[l]^(.3){\epsilon_{2k-2}} & \cdots \ar[l] & \pi_{k+ 1} \ar[l] & 
}
\]
For each $k > 0$ and a $2k$-loop $(\pi, \epsilon)$, for every $1 \leq
i < j \leq 2k$, we define
\[
\pi_{[i, j]} = \pi_i \stackrel{\epsilon_i}{\rightarrow}\pi_{i+1}
\rightarrow \cdots \stackrel{\epsilon_{j-1}}{\rightarrow} \pi_j, \,
\]
and its reverse
$\widetilde{\pi_{[i, j]} }$ to be the path
\[ 
\ \widetilde{\pi_{[i, j]}} = \pi_j
\stackrel{\epsilon_{j-1}}{\rightarrow}\pi_{j-1} \rightarrow \cdots
\stackrel{\epsilon_{i}}{\rightarrow} \pi_i\, ; 
\] 
and if $0 \leq r < s \leq 2l $ is such that $\lambda_r = \pi_j$ for
some $2l$-loop $(\lambda, \eta)$, then we define the concatenation
$\pi_{[i, j]} \circ \lambda_{[r, s]}$ to be the path
\[
\pi_{[i, j]} \circ \lambda_{[r, s]} = \pi_i
\stackrel{\epsilon_i}{\rightarrow}\pi_{i+1} \rightarrow \cdots
\stackrel{\epsilon_{j-1}}{\rightarrow} \pi_j (= \lambda_r)
\stackrel{\eta_r}{\rightarrow} \cdots
\stackrel{\eta_{s-1}}{\rightarrow} \lambda_s.
\]
Further, for any two vertices $u^+ \in \mathcal{U}^+$ and $u^- \in
\mathcal{U}^-$, we set
\[
\mathcal{E}(u^+, u^-) = \{ \epsilon \in \mathcal{E} : \epsilon\
\mbox{joins}\ u^+\ \mbox{and}\ u^-\}.
\] 

For each $k > 0$, let $P_k(\Gamma)$ be the $\mathbb{C}$-vector space
whose basis consists of loops on $\Gamma$ of length $2k$ based at
vertices in $\mathcal{U}^+$; and as for $k = 0_{\pm}$, a loop of
length $0$ is just a vertex on $\Gamma$, we set $P_{0_{\pm}}(\Gamma)$
to be the $\mathbb{C}$-vector space with $\mathcal{U}^{\pm}$ as basis.
We wish to give a planar algebra structure on the collection
$$P (\Gamma) := \{ P_k(\Gamma)\, :\, k \in\, Col\}.$$ 

\begin{defn}\label{state}
  A state $\sigma$ of a coloured tangle $T$ is a function
  \[
  \sigma : \{\mbox{regions of $T$}\} \sqcup\{\mbox{strings of $T$}\}
  \rightarrow \mathcal{U} \sqcup \mathcal{E},
  \]
  such that
  \begin{enumerate}
  \item $\sigma (\{\mbox{unshaded regions}\}) \subset \mathcal{U}^+$,
    $ \sigma (\{\mbox{shaded regions}\} ) \subset \mathcal{U}^- $;
  \item $\sigma (\{\mbox{strings} \} ) \subset \mathcal{E}$; and
  \item if a string $s$ lies in the closure of two regions $r_1$ and
    $r_2$, then $\sigma(s)$ is the edge joining the vertices $\sigma
    (r_1)$ and $\sigma (r_2)$.
  \end{enumerate}
\end{defn}

Note that a state $\sigma$ of a $k_0$-tangle $T$ with $b$ internal
boxes of colours $k_1, \ldots, k_b$, respectively, induces unique
loops at the internal as well as the external
boxes of $T$:\\
\hspace*{5mm} Suppose around the $j$-th internal box $D_j$ of $T$, say
of colour $ k_j = 4$, the state $\sigma$ looks like
\begin{center}
\psfrag{u0}{$u_0$}
\psfrag{u1}{$u_1$}
\psfrag{u2}{$u_2$}
\psfrag{u3}{$u_3$}
\psfrag{u4}{$u_4$}
\psfrag{u5}{$u_5$}
\psfrag{u6}{$u_6$}
\psfrag{u7}{$u_7$}
\psfrag{e0}{$\epsilon_0$}
\psfrag{e1}{$\epsilon_1$}
\psfrag{e2}{$\epsilon_2$}
\psfrag{e3}{$\epsilon_3$}
\psfrag{e4}{$\epsilon_4$}
\psfrag{e5}{$\epsilon_5$}
\psfrag{e6}{$\epsilon_6$}
\psfrag{e7}{$\epsilon_7$}
\psfrag{Bj}{$D_j$}
\includegraphics[scale=0.75]{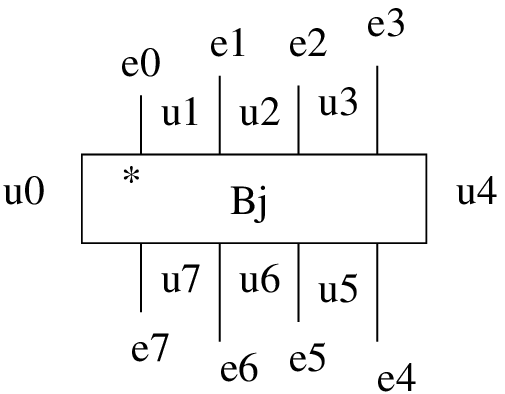},
\end{center} 
then, by condition $(3)$ of Definition \ref{state}, we get a loop
\[
\vcenter{\xymatrix @-1.2pc{ & u_1 \ar[r]^{\epsilon_1} & u_2
    \ar[r]^{\epsilon_2} & u_3 \ar[dr]^{\epsilon_3} &
    \\
    u_0 \ar[ru]^{\epsilon_0}& & & & u_4 \ar[ld]^{\epsilon_4} \\
    & u_7 \ar[lu]^{\epsilon_7} & u_6 \ar[l]^{\epsilon_6} & u_5
    \ar[l]^{\epsilon_5} & }}.
\]
And if $k_j \in \{ 0_{\pm}\}$, then the state $\sigma$ simply induces
a vertex in $\mathcal{U}^+$ or $\mathcal{U}^- $ according as $D_j$ is
a $0_+$- or $0_-$-box. We take this vertex as the $2 k_j$-loop induced
by $\sigma$ at $D_j$. Similarly, $\sigma$ also induces a $2 k_0$-loop
at the external box $D_0$.  We denote this unique $2 k_j$-loop by the
pair $(\pi^{\sigma}_{D_j}, \epsilon^{\sigma}_{D_j})$, $ 0 \leq j \leq
b$.  For each $0 \leq j \leq b $, we say that the state $\sigma$ is
compatible with a $2k_j$-loop $(\pi, \epsilon)$ at its $j$-th box
$D_j$ if $(\pi^{\sigma}_{D_j}, \epsilon^{\sigma}_{D_j}) = (\pi,
\epsilon)$. We can now define the tangle maps.

Given a $k_0$-tangle $T$ as above, we first isotope it to a ``standard
form'', i.e.,
\begin{itemize}
\item first rotate all the internal boxes of $T$ so that their
  $*$-vertices are on top left-corner; and 
\item then isotope all the strings, if necessary, so that any
  singularity of the $y$-coordinate function for strings is either a
  local maximum or a local minimum.
\end{itemize}
Given $2 k_j$-loops $\{ (\pi^{(j)}, \epsilon^{(j)}) ;\, 0
\leq j \leq b \}$ (with the convention that $2 k_j = 0_{\pm},\,
\mbox{if}\ k_j = 0_{\pm} $, and $ (\pi^{(j)}, \epsilon^{(j)})$ is an
appropriate element of $\mathcal{U}$), we define the coefficient of $
(\pi^{(0)}, \epsilon^{(0)})$ in the vector $Z^{P (\Gamma)}_T(
(\pi^{(1)}, \epsilon^{(1)}) \otimes \cdots \otimes (\pi^{(b)},
\epsilon^{(b)}))$ \footnote{We follow the convention that empty tensor
  product denotes the scalars.} by

\begin{eqnarray}\label{coefficient-1}
  \lefteqn{Z^{P (\Gamma)}_T( (\pi^{(1)}, \epsilon^{(1)}) \otimes \cdots 
    \otimes (\pi^{(b)},
    \epsilon^{(b)}))_{(\pi^{(0)},\, \epsilon^{(0)})}}\nonumber \\
  & = & \sum_{\sigma\, \in\,
    \left\{ \mbox{\begin{tabular}{c} states of $T$ \\ compatible with \\
          $(\pi^{(j)}, \epsilon^{(j)})$ at $D_j$, \\ $\forall\ 0 \leq j \leq b$
      \end{tabular} } \right\}} \prod_{\mbox{\begin{tabular}{c} singularities
        $\alpha$ \\ of $y$-coordinate \\ on strings \end{tabular}}} \mu_{\alpha} , 
  \end{eqnarray}
where $\mu_{\alpha} := \mu_x / \mu_y$ taking $x = \sigma$(inner region
at $\alpha$), $y = \sigma$(outer region at $\alpha$).
\begin{center}
  \psfrag{alpha}{$\alpha$} \psfrag{x}{$x$} \psfrag{y}{$y$}
\includegraphics[scale=0.35]{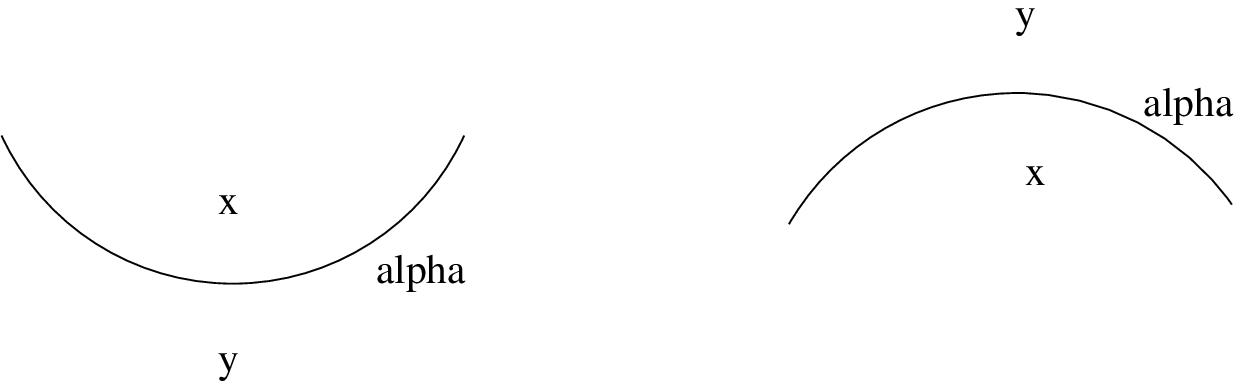}
\end{center}

This gives the required tangle map $Z^{P (\Gamma)}_T :
\otimes_{i=1}^bP_{k_i} \rightarrow P_{k_0} (\Gamma)$. The fact that
this definition does not depend upon the way we isotope the tangle $T$
to a standard form and that these tangle maps satisfy the
compatibility conditions for composition of tangles is precisely the
crux of Theorem $3.1$ in \cite{Jon00}.  \newline

\noindent {\em *-structure on $P(\Gamma)$:}

There is a natural $*$-algebra structure on each $P_k(\Gamma)$, such
that the $*$-compatibility condition (\ref{*-condition}) holds:

On $P_{0_{\pm}}(\Gamma)$, we define $(u^{\pm})^* =u^{\pm} $. And for
each $k > 0$, given a $2k$-loop $(\pi, \epsilon)$, we define
\[
\xymatrix @-1.5pc{ & & \pi_{2k-1} \ar[r]^{\epsilon_{2k-2}} &
  \pi_{2k-2} \ar[r] & \cdots \ar[r] & \pi_{k+ 1}
  \ar[rd]^{\epsilon_{k}} & \\
  (\pi, \epsilon)^* = & \pi_0 \ar[ru]^{\epsilon_{2k-1}} & & & & &
  \pi_{k} \ar[ld]^{\epsilon_{k-1}} .\\
  & & \pi_1 \ar[lu]^{\epsilon_0} & \pi_2 \ar[l]^{\epsilon_{1}} &
  \cdots \ar[l] & \pi_{k-1} \ar[l] & }
\]
This gives an involution on $P_k(\Gamma)$, with respect to the
multiplication induced by $Z^{P(\Gamma)}_{M_k}$, and the
$*$-compatibility condition (\ref{*-condition}) holds.

\subsection{Some Calculations}
For a coloured tangle $T$, it is quite illustrative to write the
coefficient $Z^{P (\Gamma)}_T( (\pi^{(1)}, \epsilon^{(1)}) \otimes
\cdots \otimes (\pi^{(b)}, \epsilon^{(b)}))_{(\pi^{(0)},\,
  \epsilon^{(0)})}$ as a picture, as shown by the self explanatory
example in equation \ref{possible-states}, where $T$ is the coloured
tangle obtained by removing the labels from the picture.

\begin{equation} \label{possible-states}
\begin{array}{rcl} Z^{P(\Gamma)}_T ((\pi, \epsilon) \otimes (\lambda,
  \eta))_{_{(\theta, \nu)}}
 & = & \vcenter{
    \mbox{
      \psfrag{p0}{$^{_{\pi_0}}$} \psfrag{p1}{$^{_{\pi_1}}$}
      \psfrag{p2}{$^{_{\pi_2}}$} \psfrag{p3}{$^{_{\pi_3}}$}
      \psfrag{p4}{$^{_{\pi_4}}$} \psfrag{p5}{$^{_{\pi_5}}$}
      \psfrag{e0}{$^{_{\epsilon_0}}$} \psfrag{e1}{$^{_{\epsilon_1}}$}
      \psfrag{e2}{$^{_{\epsilon_2}}$} \psfrag{e3}{$^{_{\epsilon_3}}$}
      \psfrag{e4}{$^{_{\epsilon_4}}$} \psfrag{e5}{$^{_{\epsilon_5}}$}
      \psfrag{l0}{$^{_{\lambda_0}}$} \psfrag{l1}{$^{_{\lambda_1}}$}
      \psfrag{l2}{$^{_{\lambda_2}}$} \psfrag{l3}{$^{_{\lambda_3}}$}
      \psfrag{f0}{$^{_{\eta_0}}$} \psfrag{f1}{$^{_{\eta_1}}$}
      \psfrag{f2}{$^{_{\eta_2}}$} \psfrag{f3}{$^{_{\eta_3}}$}
      \psfrag{v0}{$^{_{\theta_0}}$} \psfrag{v1}{$^{_{\theta_1}}$}
      \psfrag{v2}{$^{_{\theta_2}}$} \psfrag{v3}{$^{_{\theta_3}}$}
      \psfrag{v4}{$^{_{\theta_4}}$} \psfrag{v5}{$^{_{\theta_5}}$}
      \psfrag{t0}{$^{_{\nu_0}}$} \psfrag{t1}{$^{_{\nu_1}}$}
      \psfrag{t2}{$^{_{\nu_2}}$} \psfrag{t3}{$^{_{\nu_3}}$}
      \psfrag{t4}{$^{_{\nu_4}}$} \psfrag{t5}{$^{_{\nu_5}}$}
\includegraphics[scale=0.5]{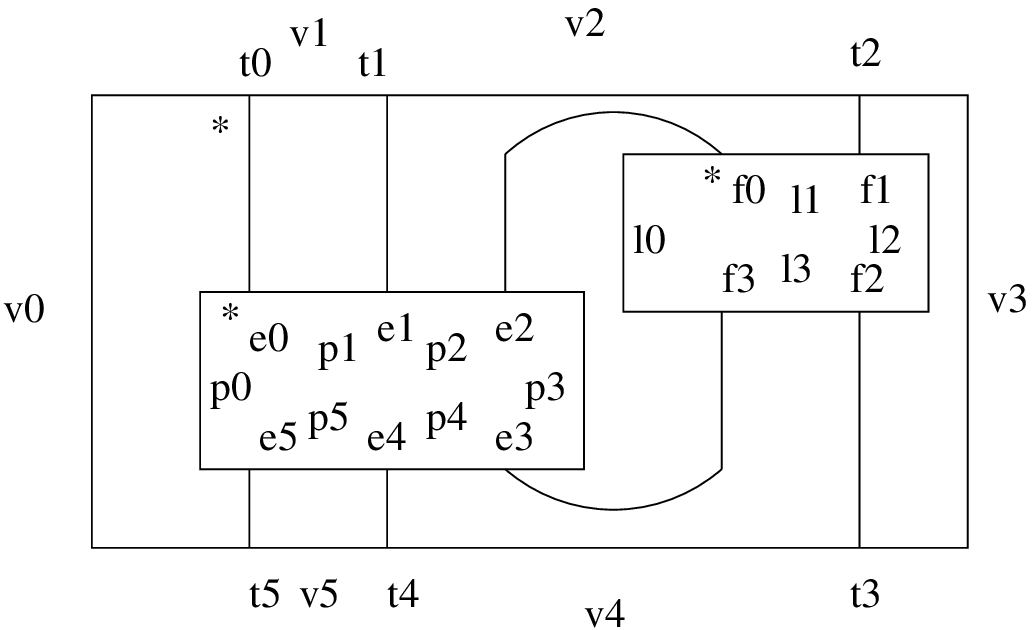}
}} \\
& = &
\delta^{\theta_{[0,\, 2]}}_{\pi_{[0,\, 2]}}\ 
\delta^{\theta_{[4,\, 6]}}_{\pi_{[4,\, 6]}}\ 
\delta^{\theta_{[2,\, 4]}}_{\lambda_{[1,\, 3]}}\ 
\delta^{\pi_{[2,\, 4]}}_{\widetilde{ \lambda_{[3,\, 4] \circ 
      \lambda_{[0,\, 1]} }}}. 
\end{array}
\end{equation}

With the help of this diagrammatic approach, we list the expressions
of some tangle maps, which we shall need in the sequel.

\begin{equation}
\begin{array}{rcl}
Z^{P(\Gamma)}_{I_k^{k+1}} ((\pi, \epsilon)) 
  & = &\sum_{\left\{
      \begin{array}{c}\mbox{\small possible}\, v \in \mathcal{U},\\
        \epsilon \in \mathcal{E}(\pi_k, v)\end{array}\right\} } 
  \vcenter{ \xymatrix @-1.5pc{ 
      & \pi_1 \ar[r]^{\epsilon_1} & \cdots
      \ar[r]^{\epsilon_{k-1}} & \pi_{k} \ar[dr]^{\epsilon} &  \\
      \pi_0 \ar[ru]^{\epsilon_0} & & & &\  v \ar[ld]^{\epsilon}\\
      & \pi_{2k-1} \ar[lu]^{\epsilon_{2k-1}} & \cdots
      \ar[l]^(.3){\epsilon_{2k-2}} & \pi_{k} \ar[l]^{\epsilon_k} &
    }},\label{I_^k+1}\\
 && \forall\ \mbox{$2k$-loops}\ (\pi, \epsilon).
\end{array}
\end{equation}
\begin{equation}
\begin{array}{rcl} 
  Z^{P(\Gamma)}_{M_{k}} ((\pi, \epsilon) \otimes (\lambda, \eta) ) & = & 
  \vcenter{  \xymatrix @-1.5pc{ 
      &  & \pi_1 \ar[r]^{\epsilon_1} & \cdots
      \ar[r] & \pi_{k-1} \ar[dr]^{\epsilon_{k-1}} &  \\
      \delta^{\widetilde{\pi_{[k,\, 2k]}}}_{\lambda_{[0, k]}} & 
      \pi_0 \ar[ru]^{\epsilon_0} & & & & \lambda_k
      \ar[ld]^{\eta_k}\\
      &  & \lambda_{2k-1} \ar[lu]^{\eta_{2k-1}} & \cdots
      \ar[l]^(.3){\eta_{2k-2}} & \lambda_{k+1} \ar[l] & 
    }}, \label{Mk}\\
  & & \forall \ \mbox{$2k$-loops}\ (\pi, \epsilon), (\lambda, \eta). 
\end{array}
\end{equation}
\begin{eqnarray}
  Z^{P(\Gamma)}_{(E')^k_{k}} ((\pi, \epsilon)) & = & 
  \delta^{\widetilde{\pi_{[2k-1,\, 2k]}}}_{\pi_{[0,\, 1]}}\,
  \frac{\mu_{\pi_{0}}^2}{\mu_{\pi_{1}}^2} 
  \sum_{\left\{\begin{array}{c}\mbox{\small possible}\, u \in \mathcal{U}^+,
        \\ \epsilon \in \mathcal{E}(u, \pi_1) \end{array}\right\}} \nonumber\\
  & & \vcenter{ \xymatrix @-1.5pc{ 
      & \pi_1 \ar[r]^{\epsilon_1} & \cdots
      \ar[r] & \pi_{k-1} \ar[dr]^{\epsilon_{k-1}} &  \\
      u \ar[ru]^{\epsilon} & & & & \pi_k
      \ar[ld]^{\epsilon_{k}}\\
      & \pi_{2k-1} \ar[lu]^{\epsilon} & \cdots
      \ar[l]^(.3){\epsilon_{2k-1}} & \pi_{k+1} \ar[l] & 
    }},\ \forall\ \mbox{$2k$-loops}\ (\pi, \epsilon).  \label{E'k}
\end{eqnarray}
\begin{eqnarray}
\qquad Z^{P(\Gamma)}_{\mathcal{E}^{k+1}}(1) & = & 
\sum_{\left\{\begin{array}{c} \mbox{\small 2(k+1)-loops}\,(\pi, \epsilon)\, :\\  
       \pi_{[0,\, k-1]} = \widetilde{\pi_{[k+3,\, 2k]}}, \\
      \pi_{[k-1,\, k]} = \widetilde{\pi_{[k,\, k+1]}},\\
      \pi_{[k+1,\, k+2]} = \widetilde{\pi_{[k+2,\, k+3]}}
      \end{array}\right\}}
  \frac{\mu_{\pi_k} \mu_{\pi_{k+2}}}{\mu_{\pi_{k-1}}^2}\, (\pi, \epsilon).
\label{Ek+1}  
\end{eqnarray}

\subsection{ Flip and Dual}

Given a bipartite graph $\Gamma =(\mathcal{U}^+, \mathcal{U}^-,
\mathcal{E}) $, we define another bipartite graph $\overline{\Gamma} =
(\mathcal{U}^-, \mathcal{U}^+, \mathcal{E})$ such that if $\epsilon$
is an edge in $\Gamma$ joining an even vertex $u_+$ with an odd vertex
$u_-$, then in $\overline{\Gamma}$, $\epsilon$ is an edge joining the
even vertex $u_-$ with the odd vertex $u_+$. We call the bipartite
graph $\overline{\Gamma}$ to be the flip of the graph $\Gamma$.

Given a finite, connected bipartite graph $\Gamma = (\mathcal{U}^+,
\mathcal{U}^-, \mathcal{E}, \mu)$, we note that its flip
$\overline{\Gamma} $ is also finite and connected, and $\mu$ is a spin
function for $\overline{\Gamma}$ as well.

\begin{thm}\label{dual-dual}
  The dual of the planar algebra of a finite, connected bipartite
  graph $\Gamma = (\mathcal{U}^+, \mathcal{U}^-, \mathcal{E}, \mu)$ is
  isomorphic to the planar algebra of its flip $\overline{\Gamma} =
  (\mathcal{U}^-, \mathcal{U}^+$, $\mathcal{E}, \mu) $, i.e.,
  $$^-P(\Gamma) \cong P (\overline{\Gamma}).$$
\end{thm}
\begin{pf}
  We have $^-P(\Gamma)_{0_{\pm}} = P_{\mp 0}(\Gamma) = \mathbb{C}
  [\mathcal{U}^{\mp}] = P_{\pm 0}(\overline{\Gamma})$; and for each $
  k > 0$,
\[
\begin{array}{ccl}
  ^-P(\Gamma)_k & = & P_k (\Gamma) = \mathbb{C} [\,\mbox{$2k$-loops
    on $\Gamma$ based at vertices in $\mathcal{U}^+$} \,],\ \mbox{and}\\
  P(\overline{\Gamma})_k & = & \mathbb{C} [\,\mbox{$2k$-loops on
    ${\Gamma}$ based at vertices in $\mathcal{U}^-$} \,].
  \end{array}
  \]
  Let $\varphi_{0_{\pm}} :\, ^-P(\Gamma)_{0_{\pm}} \rightarrow
  P_{0_{\pm }}(\overline{\Gamma})$ be the identity morphism of the
  underlying vector space; and for each $k > 0$, define $\varphi_k :\,
  ^-P(\Gamma)_{k} \rightarrow P_{k}(\overline{\Gamma})$ to be the map
  whose action on basis vectors is given by
  \begin{equation}
    \varphi_k ((\pi, \epsilon)) =  \frac{\mu_{\pi_0} \mu_{\pi_k}}{\mu_{\pi_1}
      \mu_{\pi_{k+1}}} 
    \vcenter{ \xymatrix @-1.2pc{ 
        & \pi_2 \ar[r]^{\epsilon_2} & \cdots
        \ar[r] & \pi_{k} \ar[dr]^{\epsilon_{k}} &  \\
        \pi_1 \ar[ru]^{\epsilon_1} & & & &\! \! \pi_{k+1}
        \ar[ld]^{\epsilon_{k+1}}\ ,\\
        & \pi_{0} \ar[lu]^{\epsilon_{0}} & \cdots
        \ar[l]^(.3){\epsilon_{2k-1}} & \pi_{k+ 2} \ar[l] &} }
\end{equation} 
for all $2k$-loops $(\pi, \epsilon)$ on $\Gamma$ based at vertices in
$\mathcal{U}^+$.  Clearly, $\varphi_k ((\pi, \epsilon))$ is a
$2k$-loop on ${\Gamma}$ based at a vertex in
$\mathcal{U}^-$.
 
And, as the above correspondence is a bijection between $2k$-loops on
$\Gamma$ and those on its flip, each $\varphi_k$ is a linear
isomorphism. We claim that $\varphi := \{\varphi_k;\, k\in\, Col \} $
is a planar algebra morphism from $^-P(\Gamma) $ onto $ P
(\overline{\Gamma}) $.

Let $\mathcal{T}$ be the collection of tangles $T$ for which the
equation in Definition \ref{planar-morphism} holds. By Lemma
\ref{closedness-of-tangles}, $\mathcal{T}$ is closed under composition
of tangles. Thus, in order to show that $\mathcal{T}$ contains all the
coloured tangles, it is enough to show, by Theorem \ref{gen-tangles},
that it contains the collection
  \[
  \mathcal{G}_1 = \{ 1^{0_{\pm}}\} \cup \{ I^{k+1}_{k}, M_k,
  E^k_{k+1}:\, k \in\, Col \} \cup \{ R_k : k \geq 2\}.
  \]
  The verification of this assertion is a routine, if laborious,
  exercise using the pictorial approach to tangle maps as described
  above.
\end{pf}

\subsection{Group Action on a Bipartite Graph with Spin Function}
\begin{defn}\label{planar-action}
  Let $\Gamma = (\mathcal{U}^+,\mathcal{U}^-, \mathcal{E})$ be a
  bipartite graph with a spin function $\mu$.  We say that a finite
  group $G$ acts on $\Gamma$ if
  \begin{itemize}
  \item $G$ acts on each of the sets $\mathcal{U}^+,\mathcal{U}^-$ and
    $\mathcal{E}$; such that
  \item if $\epsilon $ is an edge joining the vertices $u_+$ and
    $u_-$, then $g \epsilon$ is an edge between the vertices $g u_+$
    and $g u_-$.
  \end{itemize}
  Further, we say that $G$ acts on $(\Gamma, \mu)$ if $G$ acts on
  $\Gamma$ and the spin function $\mu$ is constant on $G$-orbits of
  $\mathcal{U}^{\pm}$.
\end{defn} 
\begin{thm}\label{gp-action-p-gamma}
  Let $\Gamma = (\mathcal{U}^+,\mathcal{U}^-, \mathcal{E})$ be a
  finite, connected bipartite graph with a spin function $\mu$, and
  suppose that a finite group $G$ acts on $(\Gamma, \mu)$. Then there
  is a canonical $G$-action on the planar algebra $P
  (\overline{\Gamma})$, and its $G$-invariant planar sub-algebra $P
  (\overline{\Gamma})^G$ is isomorphic to the dual planar algebra
  $^-(P (\Gamma)^G)$, i.e.,
  \[P (\overline{\Gamma})^G\, \cong\, ^-(P (\Gamma)^G). \]
\end{thm}
\begin{pf}
  For each $k > 0$, the action of $G$ on $(\Gamma, \mu)$ induces an
  action on the set $\{2k$-loops on $\Gamma$ based at vertices in
  $\mathcal{U}^+ \}$ given by
  \[
  g (\pi, \epsilon) = \vcenter{\xymatrix @-1.3pc{ & g \pi_1
      \ar[r]^{g\, \epsilon_1} &g \pi_2 \ar[r] & \cdots
      \ar[r] &g \pi_{k-1} \ar[rd]^{g\, \epsilon_{k-1}} &  \\
      g \pi_0 \ar[ru]^{g\, \epsilon_0} & & & &
      & g \pi_{k} \ar[ld]^{g\, \epsilon_{k}}, \\
      & g \pi_{2k-1} \ar[lu]^{g\, \epsilon_{2k-1}} & g \pi_{2k-2}
      \ar[l]^(.3){g\, \epsilon_{2k-2}} & \cdots \ar[l] &g \pi_{k+ 1}
      \ar[l] & } }
  \]
  for all $2k$-loops $(\pi, \epsilon)$ and $g \in G$.

  Already, $G$ acts on $\mathcal{U}^{\pm}$, thus we get a group action
  on $P_k (\Gamma)$, for all $k \in\, Col$. Thanks to the invariance
  of the spin function on the $G$-orbits of $\mathcal{U}^{\pm}$, Lemma
  \ref{closedness-of-tangles} and Theorem \ref{gen-tangles}, a bit of
  straight forward checking readily shows that it is in fact a
  $G$-action on the planar algebra $P (\Gamma)$.

  And the second assertion that $ P (\overline{\Gamma})^G \cong \,
  ^-(P (\Gamma)^G)$ is a mere consequence of Proposition
  \ref{G-invariant} and the fact that the morphism $\varphi =
  \{\varphi_k ;\, k \in \, Col\}$ giving the planar isomorphism $P
  (\overline{\Gamma}) \cong \, ^-P (\Gamma) $ in Theorem
  \ref{dual-dual} is in fact a $G$-map, i.e., $\varphi_k(g x) = g
  \varphi_k (x),\, \forall g \in G, x \in P_k, k \in\, Col$.
\end{pf}
\newline

We list some easily verifiable facts; which are also suggested in
\cite{Jon00}.

\begin{lem}\label{modulus-connected-irreducible}
  Let $\Gamma$ and $G$ be as in Theorem \ref{gp-action-p-gamma}. Then
  $P (\Gamma)^G$
  \begin{enumerate}
  \item is connected iff $G$ acts transitively on $\mathcal{U}^{\pm}$.
  \item has modulus $|| \Gamma ||$ if in addition to (1),
    $(\mu_{u}^2)_{u \in \mathcal{U}}$ is the Perron-Frobenius
    eigen-vector of the adjacency matrix $A_{\Gamma}$ of $\Gamma$,
    where $|| \Gamma ||$ is the norm of $A_{\Gamma}$.
  \item is irreducible iff $G$ acts transitively on the set
    $\{2$-loops on $\Gamma$ based at vertices in $\mathcal{U}^+\}$.
  \end{enumerate}
\end{lem}

\section{Planar Algebra of the Subgroup-Subfactor}
\subsection{The Bipartite Graph $\star_n$}\label{bipartite-graph}

\begin{figure}[h]
\begin{center}
  \psfrag{g1}{$^{_{H g_1}}$} \psfrag{g2}{$^{_{H g_2}}$}
  \psfrag{g3}{$^{_{H g_3}}$} \psfrag{g4}{$^{_{H g_4}}$}
  \psfrag{g5}{$^{_{H g_5}}$} \psfrag{gn}{$^{_{H g_n}}$}
  \psfrag{gn-1}{$^{_{H g_{n-1}}}$} \psfrag{gn-2}{$^{_{H g_{n-2}}}$}
  \psfrag{gn-3}{$^{_{H g_{n-3}}}$}\psfrag{*}{$\ast$}
\includegraphics[scale=0.75]{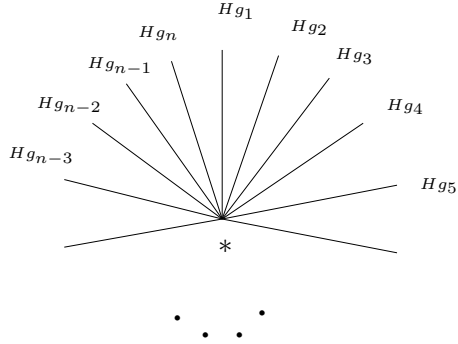}
\caption{The bipartite graph $ \star_n$}\label{n-star}
\end{center}
\end{figure}

Consider the bipartite graph $\star_n = (\mathcal{U}^+,
\mathcal{U}^-,$ $\mathcal{E})$, as in Figure \ref{n-star}, where
$\mathcal{U}^+ := \{Hg_1, \ldots, Hg_n\}$, $\mathcal{U}^- := \{*\}$
and the edge set $\mathcal{E} := \{\epsilon_1, \ldots, \epsilon_n\}$,
where $\epsilon_i$ is the edge joining the vertices $H g_i$ and $*$, $
1 \leq i \leq n$. Let $\mu : \mathcal{U}^+ \sqcup \mathcal{U}^-
\rightarrow (0, \infty)$ be the spin function given by $\mu (*) =
n^{1/4}$ and $\mu (H g_i) = 1$ for all $i \in I$.

Note that $(\mu_{u}^2)_{u\, \in\, \mathcal{U}^+\, \sqcup\,
  \mathcal{U}^-}$ is the Perron-Frobenius eigen vector of the
adjacency matrix of the bipartite graph $\star_n$, and the norm of
this adjacency matrix is $||\star_n || = \sqrt{n}$. Further, with this
set up, there is a natural action of $G$ on $(\star_n, \mu)$:

{\em $G$ acts trivially on $\mathcal{U}^-$; on $\mathcal{U}^+$ it acts
  by the natural left action, i.e., $g \cdot H g_i = H g_j \Leftrightarrow
  H g_i g^{-1} = H g_j$; and on $\mathcal{E}$ its action is induced by
  that on $\mathcal{U}^+$, i.e., $g \epsilon_i = \epsilon_j
  \Leftrightarrow g \cdot H g_i = H g_j,$ for $ g \in G, i, j \in
  I$.}
\begin{rem}\label{p-n-star-modulus}
  By Lemma \ref{modulus-connected-irreducible}, the planar algebra
  $P(\star_n)^G$ is connected and has positive modulus $\sqrt{n}$.
\end{rem}

There being only one odd vertex and no multiple edges in $\star_n$,
we can ignore the $*$-vertex and the edge joining it to an even vertex
in the notation of a loop on $\star_n$. Thus if $k$ is even, say $k
= 2r$, we simply write
\[
\left(\vcenter{\small \xymatrix @-2.5pc{
      &   H g_{i_1}, & \cdots , &\, H g_{i_{r-1}} & \\
      H g_{i_0} & & & & H g_{i_r} \\
      & H g_{i_{2r-1}}, & \cdots , &\, H g_{i_{r+1}} & }}\right)
\]
for the $2k$-loop
\[
\vcenter{\small\xymatrix @-1.1pc{ & \ast \ar [r]^{\epsilon_{i_1}} & H
    g_{i_1} \ar [r] & \cdots \ar [r] & H g_{i_{r-1}}
    \ar [r]^(.6){\epsilon_{i_{r-1}}} & \ast \ar [rd]^{\epsilon_{i_{r}}} & \\
    H g_{i_0} \ar [ru]^{\epsilon_{i_0}} & & & & & & H g_{i_r}
    \ar [dl]^{\epsilon_{i_r}},\\
    & \ast \ar [lu]^{\epsilon_{i_0}} & H g_{i_{2r-1}} \ar
    [l]^{\epsilon_{i_{2r-1}}} & \cdots \ar [l] & H g_{i_{r+ 1}} \ar
    [l] & \ast \ar[l]^(.1){\epsilon_{i_{r+1}}} & }}
\]
and if $k$ is odd, say $k =
2s+1$, we denote by
\[
\left(
\vcenter{\small\xymatrix @-2.5pc{
&   H g_{i_1}, & \cdots, &\, H g_{i_{s}}  \\
H g_{i_0} & & & \\
& H g_{i_{2s}}, & \cdots, &\, H g_{i_{s+1}}
}}
\right)
\]
the $2k$-loop 
\[\vcenter{\small \xymatrix @-1.1pc{ & \ast \ar [r]^{\epsilon_{i_1}} &
    H g_{i_1} \ar [r] & \cdots \ar [r] & \ast \ar [r]^{\epsilon_{i_s}}
    &    H g_{i_{s}} \ar [rd]^{\epsilon_{i_s}}&  \\
    H g_{i_0} \ar [ur]^{\epsilon_{i_0}}& & & & & & \ast
    \ar [ld]^{\epsilon_{s+1}},\\
    & \ast \ar [lu]^{\epsilon_{i_0}} & H g_{i_{2s}}\ar
    [l]^{\epsilon_{i_{2s}}} & \cdots \ar [l] &\ast \ar [l] & H
    g_{i_{s+1}} \ar [l]^{\epsilon_{s+1}} & }}
\]
for all $\underline{i} = (i_0, i_1, \ldots, i_{k-1}) \in I^{k}$.

We shall also view these new descriptions of $2k$-loops as elements of
$(H \backslash G)^{k}$. Also note that the action $\beta^1$ of $G$ on
$I$ is basically the natural left action of $G$ on the set $H
\backslash G = \{ Hg_1, \ldots, H g_n \}$ \footnote{Our running
  notation for orbit representatives is not be confused with $H
  \backslash G$.}.

With these simplified notations we have the following set of bases.

\begin{lem}\label{bases-p-n-star} $\{\sum_{i \in I} H g_i \}$ and 
  $\{\ast\}$ form bases for $P_{0_+}(\star_n)^G$ and\\ $P_{0_-}(
  \star_n)^G$, respectively.  For every choice of orbit representatives
  $G \backslash (H \backslash G)^{k}$ under diagonal $\beta^1$-action,
  \begin{eqnarray*} \lefteqn{ \left\{ \sum_{g \in G} g \left(
          \vcenter{\small\xymatrix @-2.5pc{
              & H g_{i_1}, & \cdots, &\, H g_{i_{r-1}} & \\
              H g_{i_0} & & & & H g_{i_{r}}\\
              & H g_{i_{2r-1}}, & \cdots, &\, H g_{i_{r+1}} & }}
        \right)
      \right.: } \\
    &\qquad \qquad \left. \left(\vcenter{\small\xymatrix @-2.5pc{
            & H g_{i_1}, & \cdots, &\, H g_{i_{r-1}} & \\
            H g_{i_0} & & &  & H g_{i_{r}}\\
            & H g_{i_{2r-1}}, & \cdots, &\, H g_{i_{r+1}} & }} \right)
      \in G \backslash (H \backslash G)^{2r} \right\} &
\end{eqnarray*}
forms a basis for $P_{2r} (\star_n)^G$ for all $r \geq 1$; and
\[\!
 \left\{ \sum_{g \in G}g  \left(\!
        \vcenter{\small\xymatrix @-2.5pc{
            & H g_{i_1}, & \cdots, &\, H g_{i_{r}}  \\
            H g_{i_0} & &  \\
            & H g_{i_{2r}}, & \cdots, &\, H g_{i_{r+1}}
          }}\! \right)\!: \!\left(\! \vcenter{\small\xymatrix @-2.5pc{
          & H g_{i_1}, & \cdots, &\, H g_{i_{r}}  \\
          H g_{i_0} & &   \\
          & H g_{i_{2r}}, & \cdots, &\, H g_{i_{r+1}}}}\! \right) \in G
    \backslash (H \backslash G)^{2r+1}\!\! \right\} 
\]
forms a basis for $P_{2r+1} (\star_n)^G$ for all $r \geq 0$.
\end{lem}

Remark \ref{p-n-star-modulus}, Lemma \ref{bases-p-n-star}, Lemma
\ref{irreducible-spherical}, and Remark
\ref{pictorial-trace-conditional} give the following:

\begin{lem}\label{p-n-star-spherical}
  With notations as in $\S \ref{bipartite-graph}$, the planar algebra
  $P(\star_n)^G$ is irreducible, spherical and admits a global trace.
\end{lem}

Moving towards the desired identification, we define
$\varphi_{0_{\pm}} : P^{R \rtimes H \subset R \rtimes G}_{0_{\pm}}
\rightarrow P_{0_{\pm}}(\star_n)^G $ by
\begin{equation}\label{varphi-0-pm}
  P^{R \rtimes H
    \subset R \rtimes G}_{0_{\pm}} := \mathbb{C} \ni\ \lambda \
  \stackrel{\varphi_{0_{\pm}}}{\longmapsto}\ \lambda 1_{0_{\pm}}\ \in
  P_{0_{\pm}}(\star_n)^G,
\end{equation}
where $1_{0_+} = \sum_{i \in I} H g_i$ and $1_{0_-} = \ast$ are the
respective multiplicative units of $P_{0_{\pm}}(\star_n)^G$; and
$\varphi_{1} : P^{R \rtimes H \subset R \rtimes G}_{1} = \mathbb{C}
\rightarrow P_{1}(\star_n)^G $ by
\begin{equation}\label{varphi-1}
  P^{R \rtimes H
    \subset R \rtimes G}_{1} = \mathbb{C} \ni\ \lambda \
  \stackrel{\varphi_{1}}{\longmapsto}\ \lambda 1_{1}\ \in  P_{1}(\star_n)^G,
\end{equation}
where $1_1$ is the multiplicative unit $\sum_{i \in I} ( H g_i )$ of $
P_{1}(\star_n)^G$.

For $k \geq 2$, we define $\varphi_{k} : P^{R \rtimes H \subset R
  \rtimes G}_{k} \rightarrow P_{k}(\star_n)^G$ on the basis vectors
described in Lemma \ref{rel-com-bases} by

\begin{equation}\label{varphi-k-odd}
\begin{array}{l}
  \lefteqn{ P^{R \rtimes H \subset R \rtimes G}_{2k+1} = N'\cap M_{2k}
    \ni \sum_{h \in H}
    h [\underline{i},\, \underline{j}]^{ev}} \\ 
  \qquad  \stackrel{\varphi_{2k+1}}{\longmapsto} 
  \sum_{g \in G} g \left(
    \vcenter{\small\xymatrix @-2.5pc{
        & H g_{i_k}, &\, H g_{i_{k-1}} g_{i_k}, &\, \cdots,&\,
        H \sqcap g_{\underline{i}}  \\  
        H & & & &  \\
        & H g_{j_k}, &\, H g_{j_{k-1}} g_{j_k}, &\, \cdots,&\, 
        H \sqcap g_{\underline{j}}    
      }}
  \right) \in P_{2k+1}(\star_n)^G,   
\end{array}
\end{equation}
for all $ (\underline{i},\, \underline{j}) \in H \backslash (I^k
\times I^k)$, $ k \geq 1$; and
\begin{equation}\label{varphi-k-even}
\begin{array}{l}
  \lefteqn{ P^{R \rtimes H \subset R \rtimes G}_{2k} = N'\cap M_{2k-1} \ni 
    \sum_{h \in H}  h [\underline{i},\, \underline{j}]^{od}} \\ 
  \qquad \stackrel{\varphi_{2k}}{\longmapsto} 
  \sum_{g \in G} g \left(
    \vcenter{\small\xymatrix @-2.5pc{
        & H g_{i_k}, &\, \cdots, &\, H  g_{i_2} \cdots g_{i_k} &  \\  
        H & & & &  H \sqcap g_{\underline{i}} \\
        & H g_{j_k}, &\, \cdots, &\, H g_{j_2} \cdots g_{j_k} &
      }}
  \right) \in P_{2k}(\star_n)^G,   
\end{array}
\end{equation}
for all $ (\underline{i},\, \underline{j}) \in H \backslash Y_k$, $
k \geq 1$.

\vspace*{3mm}{\em We want to show that
\[
\{\varphi_k : k \in Col\} =: \varphi : P^{R \rtimes H \subset R
  \rtimes G} \rightarrow P(\star_n)^G
\]
is a planar algebra isomorphism.}

\subsection{Ingredients of  the Proof}

\begin{lem}\label{dimensions} With running notations, we have
 \[ 
 dim\, P^{R \rtimes H \subset R \rtimes G}_k = dim\, P_k(\star_n)^G,\
 \forall\ k \in Col.
\]
\end{lem}
\begin{pf}
  We have already seen that $P_{0_{\pm}}(\star_n)^G$,
  $P_{1}(\star_n)^G$, $P^{R \rtimes H \subset R \rtimes G}_{0_{\pm}}$
  and $P^{R \rtimes H \subset R \rtimes G}_{1} $ are all one
  dimensional.

  For each $k \geq 1$, by Lemma \ref{rel-com-bases}, $dim\, N' \cap
  M_{2k} $ (resp., $dim\, N' \cap M_{2k-1}$) is same as the number of
  $H$-orbits of $I^k \times I^k$ (resp., $Y_k$) under the diagonal
  $\beta^k$-action.  On the other hand, by Lemma \ref{bases-p-n-star},
  $dim\, P_{2k+1}(\star_n)^G$ (resp., $dim\, P_{2k}(\star_n)^G$) is
  equal to the number of $G$-orbits of $ (H \backslash G)^{2k+1} $
  (resp., $ (H \backslash G)^{2k}$) under the diagonal
  $\beta^1$-action.

  The corresponding dimensions are same because the number of
  $H$-orbits of $I^k \times I^k$ (resp., $Y_k$) under the diagonal
  $\beta^k$-action is equal to the number of $G$-orbits of $ (H
  \backslash G)^{2k+1}$ (resp., $(H \backslash G)^{2k}$) under the
  diagonal $\beta^1$-actions.  We include a proof of the assertion in
  parentheses, and the other follows similarly.

  Note that the number of $G$-orbits of $(H \backslash G)^{2k}$ is 
  same as the number of $H$-orbits of $(H \backslash G)^{2k-1}$ under
  the diagonal actions mentioned in the previous paragraph.  This is
  true because any element of $(H \backslash G)^{2k}$ lies in the
  $G$-orbit of an element of the type $(H, H g_{i_1}, \ldots, H
  g_{i_{2k-1}})$; and
\[ 
g (H, H g_{i_1}, \ldots, H g_{i_{2k-1}}) = (H, H g_{j_1}, \ldots, H
g_{j_{2k-1}})
\]
 if and only if $ g \in H $ and
  $g ( H g_{i_1}, \ldots, H g_{i_{2k-1}}) = ( H g_{j_1}, \ldots, H
  g_{j_{2k-1}})$. And thus the assertion in the parentheses follows
  from the facts that the correspondence
  \[
  Y_k \ni (\underline{i},\,\underline{j} ) \mapsto \left(\!
    \begin{array}{cl} H g_{i_k},\, H g_{i_{k-1}} g_{i_k},\, \ldots,\,
      H g_{i_2}
      \cdots g_{i_k}, & \\
      &\!\! H \sqcap g_{\underline{i}}\\
      H g_{j_k},\, H g_{j_{k-1}} g_{j_k},\, \ldots,\, H g_{j_2} \cdots
      g_{j_k}, &
\end{array}\!\! \right) \in 
  (H \backslash G)^{2k-1},
  \]
  by the definition of the action $\beta^k$, is an $H$-injection, and
  that $|Y_k| = |X^{2k-1}|$.
\end{pf}
\newline

In particular, this gives:
\begin{lem}\label{varphi-isomorphism}
  With notations as in equations
  (\ref{varphi-0-pm}$-$\ref{varphi-k-even}), for each $k \in Col$, the
  map $\varphi_k : P_k^{R \rtimes H \subset R \rtimes G} \rightarrow
  P_k (\star_n)^G$ is a linear isomorphism.
\end{lem}
We conclude this subsection with two lemmas that we shall need for
verifying that the $\varphi_k$'s are equivariant with respect to the
inclusion and conditional expectation tangle maps.

In order to understand how the inclusion tangles act on the
subgroup-subfactor planar algebra, and the way we defined the maps
$\varphi_k$, we need to analyse the inclusions in terms of the bases
that we obtained in Lemma \ref{rel-com-bases}.  For $k$ odd, say $k =
2r-1$, it is the usual matrix algebra inclusion $\theta^{(r)}(N)' \cap
M_{I^r}(N) \subset \theta^{(r)}(N)' \cap M_{I^r}(M)$, whereas for $k$
even, the inclusion is given by an appropriate $\Theta_r$.

\begin{lem}\label{inclusion-even-odd}
  For each $k \geq 1$, the inclusion $N' \cap M_{2k} \subset N'\cap
  M_{2k+1}$ is given as under:
\begin{equation}
\begin{array}{l}
  \theta^{(k)}(N)' \cap M_{I^k}(M) \ni \sum_{h \in\, H} h [\underline{i},\, 
  \underline{j}]^{ev}  \\
  \qquad  \stackrel{\Theta_{k+1}}{\longmapsto}   
  \sum_{\left\{ h \in\, H,\, r,\, s\, \in\, I\ : \atop{((r,\underline{i}),\,
        (s,\underline{j}))\, 
        \in\, Y_{k+1}} \right\}} 
  h [(r,\underline{i}),\, (s,\underline{j})]^{od}\, \in\ \theta^{(k+1)}(N)'
  \cap M_{I^{k+1}}(N),
\end{array}
\end{equation}
for all $(\underline{i},\, \underline{j}) \in H \backslash (I^k \times
I^k)$.
\end{lem}
\begin{pf}
  First we note that, for each $k \geq 1$ and $(\underline{i},\,
  \underline{j}) \in I^k \times I^k$,
  \[
\Theta_{k+1} ([\underline{i},\, \underline{j}]^{ev}) = \sum_{r, s\,
    \in I: \atop{((r,\underline{i}),\, (s,\underline{j})) \in\,
      Y_{k+1}}} [(r,\underline{i}),\, (s,\underline{j})]^{od} .
\]
Indeed, by the description of $\Theta_{k+1}$ as in Lemma
\ref{Theta_k}, we have
\[
\begin{array}{rcl} \Theta_{k+1} ([\underline{i},\,
  \underline{j}]^{ev})_{\underline{u},\, \underline{v}} & = &
  \theta_{u_1,\, v_1} ([\underline{i},\,
  \underline{j}]^{ev}_{\underline{u}_{[2},\, \underline{v}_{[2} } ) \\
  & = & \theta_{u_1,\, v_1}
  (\delta^{\underline{i}}_{\underline{u}_{[2}}\,
  \delta^{\underline{j}}_{\underline{v}_{[2}}\, u_{(\sqcap
    g_{\underline{i}})(\sqcap g_{\underline{j}})^{-1}}) \\
  & = & \delta^{\underline{i}}_{\underline{u}_{[2}}\,
  \delta^{\underline{j}}_{\underline{v}_{[2}}\, 1_{H} (g_{u_1}(\sqcap
  g_{\underline{i}})(\sqcap g_{\underline{j}})^{-1} g_{v_1}^{-1})\
  u_{(\sqcap g_{(u_1, \underline{i})})(\sqcap g_{(v_1,
      \underline{j})})^{-1}},
 \end{array}
\]
$ \forall\, \underline{u},\, \underline{v} \in I^{k+1}$. Thus we
just need to show that
\[
  \Theta_{k+1} (h [\underline{i},\, \underline{j}]^{ev})  =  \sum_{ r,\,
    s\, \in\, I\ : \atop{((r,\underline{i}),\, (s,\underline{j}))\, \in\,
      Y_{k+1}}} h [(r,\underline{i}),\, (s,\underline{j})]^{od},\
    \forall\, h \in H,
\] 
\begin{eqnarray}
  \lefteqn{ \mbox{i.e.,}\ \sum_{x,\, y\,
      \in\, I\ : \atop{((x, \beta^{k}_h (\underline{i})),\, (y,
        \beta^{k}_h (\underline{j})))\, \in\, Y_{k+1}}} [(x,
    \beta^{k}_h (\underline{i})),\,
    (y, \beta^{k}_h (\underline{j}) )]^{od}}\nonumber \\
  & = & \sum_{ r,\, s\, \in\, I\ :
    \atop{((r,\underline{i}),\, (s,\underline{j}))\, \in\, Y_{k+1}}}
  [\beta^{k+1}_h (r,\underline{i}),\, \beta^{k+1}_h
  (s,\underline{j})]^{od},\ \forall\, h \in H. \label{eqn1}
\end{eqnarray}

For each $h \in H$ and $ (\underline{u},\, \underline{v}) \in Y_{k+1}
$, the coefficient of $[\underline{u},\, \underline{v}]^{od}$ on
$L.H.S.$ (resp., $R.H.S.$) of equation (\ref{eqn1}) is the number of
elements in the set
\[
\begin{array}{l}
  L_{\underline{u},\, \underline{v}} := \{ (x, y) \in I\times I : ((x,
  \beta^{k}_h (\underline{i})),\, (y, \beta^{k}_h (\underline{j})))
  \in\, Y_{k+1},\\
  \qquad \qquad \qquad  { [(x, \beta^{k}_h (\underline{i})),\, 
    (y, \beta^{k}_h (\underline{j}))]^{od} }
  = [\underline{u},\, \underline{v}]^{od} \}
\end{array}
\]
\[
\begin{array}{l}\mbox{(resp.,}\ R_{\underline{u},\, \underline{v}} :=
  \{ (r, s) \in I \times I : r,\, s\, \in\, I\ : ((r,\underline{i}),\,
  (s,\underline{j})) \in\, Y_{k+1},\\
  \qquad \qquad \qquad \qquad {[\beta^{k+1}_h ((r,\underline{i})),\,
    \beta^{k+1}_h ((s,\underline{j}))]^{od} = [\underline{u},\,
    \underline{v}]^{od} \}}\mbox{)}.
\end{array}
\]
Viewing the two sides of equation (\ref{eqn1}) as elements of $R'\cap
M_{2k+1}$, we shall be done once we show that these coefficients are
same.

Note that if $ ( \beta^{k+1}_h ((r,\underline{i})), \, \beta^{k+1}_h
((s,\underline{j}) )) = (\underline{u},\, \underline{v}) $, for some
$h \in H$ and $( (r,\underline{i}),\, (s,\underline{j}) ) \in
Y_{k+1}$, then $\underline{u}_{[2} = \beta^{k}_h (\underline{i})$,
$\underline{v}_{[2} = \beta^{k}_h (\underline{j})$ and
\[
((u_1,\beta^{k}_h (\underline{i}) ),\, (v_1, \beta^{k}_h
(\underline{j}))) \in Y_{k+1}.
\]
Consider the maps $ \vcenter{\small\xymatrix{ L_{\underline{u},\,
      \underline{v}} \ar @<0.5ex> [r]^{\phi} & R_{\underline{u},\,
      \underline{v}} \ar @<0.5ex> [l]^{\psi}}} $ given as under:

For each $(x, y) \in L_{\underline{u},\, \underline{v}}$ (resp., $(r,
s) \in R_{\underline{u},\, \underline{v}}$ ), $\phi ((x, y )) = (p, q)
\in R_{\underline{u},\, \underline{v}}$ (resp., $\psi ((r, s)) = (w,
z)$ ), where $(r, s)$ (resp., $(w, z)$) is given by the equations
\[
\left\{ \begin{array}{l} H g_p = H g_x (\sqcap g_{\beta^{k}_h
      (\underline{i})}) h (\sqcap g_{\underline{i}})^{-1}\
    \mbox{and} \\
    H g_q = H g_y (\sqcap g_{\beta^{k}_h (\underline{j})}) h (\sqcap
    g_{\underline{j}})^{-1}
  \end{array}\right.
\]
\[
\mbox{(resp.,}\ \beta^{k+1}_h ((r,\underline{i})) = (w, \beta^k_h
(\underline{i})) \ \mbox{and}\ \beta^{k+1}_h ((s,\underline{j})) = (z,
\beta^k_h (\underline{j})) \mbox{)}.
 \]
 It can be seen that the maps $\phi$ and $\psi $ are inverses of each
 other. Thus the coefficient of $[\underline{u},\,
 \underline{v}]^{od}$ is same on both sides of equation (\ref{eqn1}),
 and we are done.
\end{pf}
\newline

The following lemma is a consequence of the uniqueness of trace
preserving conditional expectation.
\begin{lem}\label{conditional-expectation-fact}
  Let $A_0 \subset A_1 $ and $ B_0 \subset B_1 $ be inclusions of
  finite dimensional $C^*$-algebras. Suppose $tr_{A_1}$ and $tr_{B_1}$
  are faithful traces on $A_1$ and $B_1$, respectively. Let $\phi_{i}
  : A_i \rightarrow B_i,\ i = 0, 1$ be $C^*$-isomorphisms preserving
  the traces such that $\phi_1 \slash A_0 = \phi_0$. Then
  $E^{B_1}_{B_0} \circ \phi_1 = \phi_0 \circ E^{A_1}_{A_0}$, where
  $E^{X_1}_{X_0} : X_1 \rightarrow X_0$ is the unique $tr_{X_1}$-trace
  preserving conditional expectation for $X_i \in \{A_i, B_i \}, i =
  0, 1$.
\end{lem}

\subsection{The Identification}
We need to show that $\varphi$, as defined in equations
(\ref{varphi-0-pm}$-$\ref{varphi-k-even}), commutes with all the
tangle maps. 

{\em For notational convenience, in what follows, we shall simply
  write $P$ (resp., $Q$) for the subgroup-subfactor planar algebra
  $P^{R \rtimes H \subset R \rtimes G}$ (resp., the invariant planar
  subalgebra $P(\star_n)^G$). }

\begin{lem}\label{inclusion-lemma}
  $ \varphi_{k+1} \circ Z^P_{I_{k}^{k+1}} = Z^{Q
  }_{I_{k}^{k+1}}\circ \varphi_k ,\ \forall \, k \in Col.  $
\end{lem}
\begin{pf}
  We clearly have $ \varphi_1 \circ Z^{P}_{I_{0_{\pm}}^1} = Z^{
    Q}_{I_{0_{\pm}}^1} \circ \varphi_{0_{\pm}}$, by
  (\ref{varphi-0-pm}) and (\ref{varphi-1}).

  Now for $k \geq 1$, we give a proof for odd $k$, the even case is
  similar. Suppose $k = 2r+1$ for some $ r \geq 0$.  Then, by Lemma
  \ref{inclusion-even-odd}, $Z^{P}_{I_{2r+1}^{2r+2}} : N'\cap M_{2r}
  \hookrightarrow N'\cap M_{2r+1}$ is given by the inclusion map
  $\Theta_{r+1}$. On the other hand, by (\ref{I_^k+1}) for the
  bipartite graph $\star_n$, we have
\[
\begin{array}{l}
  Q_{2r+1}  \ni \sum_{g\, \in\, G} g \left( 
    \vcenter{\small \xymatrix  @-2.5pc{
        & H g_{i_1}, & \cdots, &\, H g_{i_{r}}  \\
        H g_{i_0} & & &  \\
        & H g_{i_{2r}}, & \cdots, &\, H g_{i_{r+1}} }} \right)    \\
  \  \stackrel{Z^{Q}_{I_{2r+1}^{2r+2}}}{\longmapsto} 
  \sum_{j\, \in\, I,\,  g\, \in\, G}
  g \left( \vcenter{\small \xymatrix @-2.5pc{
        & H g_{i_1}, & \cdots, &\, H g_{i_r} &  \\
        H g_{i_0} & & & & H g_{j} \\
        & H g_{i_{2r-1}}, & \cdots, &\, H g_{i_{r+1}} & }}\right)\ \in\, 
  Q_{2r+2}.
\end{array}
\]

If $r = 0$, then for each $\lambda \in \mathbb{C} = P_1 $, by
(\ref{varphi-1}) and (\ref{varphi-k-even}), we have

\begin{eqnarray*}
  (Z^{Q}_{I_{1}^2} \circ \varphi_{1})(\lambda) & = & 
  Z^{Q}_{I_{1}^2} (\lambda \sum_{i \in I}  (Hg_i))\\
  & = & Z^{Q}_{I_{1}^2} \left( \frac{\lambda}{|\mbox{Iso}_H(H)|} 
    \sum_{g\, \in\, G} g (H)\right) \\
  & = &
  \frac{\lambda}{|H|} \sum_{i\,
    \in\, I,\, g\, \in\, G} g (H,\, H g_i ) 
\end{eqnarray*}  and
\begin{eqnarray*}
(\varphi_2 \circ Z^P_{I_{1}^2}) (\lambda)&  = & 
  \varphi_2 \left(\frac{\lambda}{| \mbox{Iso}_H((1,1)) |}
    \sum_{g\, \in\, G}  g [1, 1]^{od}\right)\\
  &  = & \frac{\lambda}{|H|}\, \varphi_2 \left(\sum_{i\, \in\, I,\, h\, \in\, H} 
    h g_i [1,\, 1]^{od}\right)\\
  & = &\frac{ \lambda}{|H|} \sum_{i\, \in\, I,\, g\, \in\, G} g (H,\, H g_i ),
\end{eqnarray*}
showing that $ \varphi_2 \circ Z^{P}_{I_1^2} = Z^Q_{I_{1}^1} \circ
\varphi_{1}$, where for an $H$-set $X$ and $x \in X,\ \mbox{Iso}_H(x)
:= \{h \in H: hx = x \}$.

And for $r \geq 1$, for each $(\underline{i},\, \underline{j}) \in H
\backslash (I^r \times I^r) $, by Lemma \ref{inclusion-even-odd}, we
have

\begin{eqnarray*}
\lefteqn{  (\varphi_{2r+2} \circ Z^{P}_{I_{2r+1}^{2r+2}}) 
    \left(\sum_{h\, \in\, H} h
      [\underline{i},\, \underline{j}]^{ev}\right)}\\
  &= &\varphi_{2r+2} \left(\sum_{ h\, \in\,
      H,\, x, y \, \in\, I: \atop{((x,\,\underline{i}),\, 
        (y,\, \underline{j}))\, \in\, Y_{r+1}}}
    h [(x, \underline{i}),\, (y, \underline{j})]^{od}\right) \\
  &=& \sum_{g\, \in\, G,\,  x\, \in\, I} 
  g  \left( \vcenter{\small\small \xymatrix @-2.5pc{
        & H g_{i_r}, & \cdots, &\, H g_{i_1} \cdots g_{i_r} &  \\
        H  & & & & H \sqcap g_{(x,\, \underline{i})}  \\
        & H g_{j_r}, & \cdots, &\,  H g_{j_1} \cdots g_{j_r} & 
      }}\right)\\
\mbox{and}\quad  \lefteqn{ (Z^{Q}_{I_{2r+1}^{2r+2}} \circ \varphi_{2r+1} )
    \left(\sum_{h \in H} h
      [\underline{i},\, \underline{j}]^{ev}\right)} \\
  &=& Z^{Q}_{I_{2r+1}^{2r+2}} 
  \left( \sum_{g \in G} g  \left( \vcenter{\small \xymatrix @-2.5pc{
          & H g_{i_k}, & \cdots,  &\,    H \sqcap g_{\underline{i}} \\
          H  &  & & \\
          & H g_{j_k}, & \cdots, &\,   H \sqcap g_{\underline{j}}
        }}\right) \right)\\ 
  &= &\sum_{x \in I,\, g \in G} g  \left( \vcenter{\small \xymatrix @-2.5pc{
        & H g_{i_k}, & \cdots,  &\,  H \sqcap g_{\underline{i}} &  \\
        H  & & &  & H g_x\\
        & H g_{j_k}, & \cdots,  &\,   H \sqcap g_{\underline{j}} & 
      }}\right) \\
& =& \sum_{j \in I,\, g \in G} g  \left( \vcenter{\small \xymatrix @-2.5pc{
        & H g_{i_k}, & \cdots, &\, H g_{i_1} \cdots g_{i_k} &  \\
        H  & & & & H \sqcap g_{(j,\underline{i})}  \\
        & H g_{j_k}, & \cdots, &\,  H g_{j_1} \cdots g_{j_k} & 
      }}\right).
\end{eqnarray*}
Thus $\varphi_{k+1} \circ Z^P_{I_{k}^{k+1}} = Z^{Q}_{I_{k}^{k+1}}\circ
\varphi_k $.
\end{pf}

\begin{lem}\label{multiplication-lemma}
  $ \varphi_{k} \circ Z^P_{M_{k}} = Z^{Q}_{M_{k}}\circ (\varphi_k
  \otimes \varphi_k),\ \forall\, k \in Col.  $
\end{lem}
\begin{pf}
  There is nothing to prove for $k = 0_{\pm}$. Both $P$ and $Q$ being
  irreducible, there is nothing to be proved for $k = 1$ either. For
  $k \geq 2$, we give a proof for the case when $k$ is odd, say $k =
  2r + 1$, for some $r \geq 1$, and the other case can be proved on
  exactly similar lines. Let
\[
X_{2r+1} =\left\{ \left(\vcenter{\small \xymatrix @-2.5pc{
        &  H g_{i_1}, &\, \cdots, &\,  H g_{i_r} \\
        H  & & & \\
        & H g_{j_1}, &\, \cdots, &\, H g_{j_r} }} \right) \in ( H
  \slash G)^{2r+1} : \underline{i}, \underline{j} \in I^r
\right \} .
\] 

Then $X_{2r+1}$ is invariant under the diagonal $\beta^1$-action of
$H$ on $(H \backslash G)^{2r+1}$ and it can also be identified with
$(H \backslash G)^{2r}$ as $H$-sets.  Thus, as in the proof of Lemma
\ref{dimensions}, the correspondence
\[
I^r \times I^r\, \ni\, ( \underline{i},\, \underline{j}) \longmapsto
\left(\vcenter{\small \xymatrix @-2.5pc{ & H g_{i_r},&\, Hg_{i_{r-1}}
      g_{i_r}, &\, \cdots, &\, H \sqcap g_{\underline{i}} \\
      H & & & \\
      & H g_{j_r}, &\, H g_{j_{r-1}}g_{j_r}, &\, \cdots, &\, H \sqcap
      g_{\underline{j}} }} \right)\, \in\, X_{2r+1}
\] 
is an $H$-bijection; and a set of representatives of $H$-orbits of
$X_{2r+1}$ is also a set of representatives of $G$-orbits of $(H
\backslash G )^{2r+1}$.  This shows that $ \left\{
  \left(\vcenter{\small \xymatrix @-2.5pc{ & H g_{i_r},
        &\,Hg_{i_{r-1}} g_{i_r}, &\, \cdots, &\,
        H \sqcap g_{\underline{i}} \\
        H & & & \\
        & H g_{j_r}, &\, Hg_{j_{r-1}} g_{j_r}, &\, \cdots, &\, H
        \sqcap g_{\underline{j}} }} \right) : ( \underline{i},\,
  \underline{j}) \in H \backslash ( I^r \times I^r) \right\}$ is a set
of representatives of $G$-orbits of $(H \backslash G)^{2r+1}$. In
particular, by Lemma \ref{bases-p-n-star},\\ $\left\{ \sum_{g \in G}g
  \left(\vcenter{\small \xymatrix @-2.4pc{ & H g_{i_r},
        &\,Hg_{i_{r-1}} g_{i_r}, &\, \cdots, &\,
        H \sqcap g_{\underline{i}} \\
        H & & & \\
        & H g_{j_r}, &\, Hg_{j_{r-1}} g_{j_r}, &\, \cdots, &\, H
        \sqcap g_{\underline{j}} }} \right) : ( \underline{i},\,
  \underline{j}) \in H \backslash ( I^r \times I^r) \right\}$ forms a
basis for $Q_{2r+1} $.

We now prove that $\varphi$ commutes with the tangle map for the
multiplication tangle $M_{2r+1}$. For each pair of representatives
$(\underline{i},\, \underline{j}),\, (\underline{\tilde{i}},\,
\underline{\tilde{j}}) \in H \backslash (I^r \times I^r)$, we have
\begin{eqnarray*}
  \lefteqn{ \left(\sum_{h\, \in\, H} h [\underline{i},\,
      \underline{j}]^{ev}\right) 
    \left(\sum_{\tilde{h}\, \in\, H} \tilde{h} [\underline{\tilde{i}},\,
      \underline{\tilde{j}}]^{ev}\right)
    =  \sum_{h,\, \tilde{h}\, \in\, H}
    \delta^{\beta^r_h (\underline{j})}_{\beta^r_{\tilde{ h}}
      (\underline{\tilde{i}})}\ 
    [\beta^r_h (\underline{i}),\,
    \beta^r_{\tilde{ h}} (\underline{\tilde{j}})]^{ev} }\\
  & = & \sum_{\hat{h}\, \in\, H,\, (\underline{\hat{i}},\, 
    \underline{\hat{j}}) \in\, H \backslash    (I^r \times I^r)}   
  C_{\underline{\hat{i}},\, \underline{\hat{j}}}\ \hat{h} [\underline{\hat{i}},\, 
  \underline{\hat{j}}]^{ev}\\
  &  \stackrel{\varphi_{2r+1}}{\longmapsto}  & \sum_{\hat{g}\,
    \in\, G,\,(\underline{\hat{i}},\, 
    \underline{\hat{j}})\, \in\, H \backslash  (I^r \times I^r)}  
  C_{\underline{\hat{i}},\, \underline{\hat{j}}}\ \hat{g} 
  \left( \vcenter{\small\xymatrix @-2.5pc{  
        & H g_{\hat{i}_r}, &\, H g_{\hat{i}_{r-1} } g_{\hat{i}_r}, &\, \cdots, 
        &\, H \sqcap g_{\underline{\hat{i}}} \\
        H & & & \\
        & H g_{\hat{j}_r}, &\, Hg_{\hat{j}_{r-1}} g_{\hat{j}_r}, &\, \cdots, 
        &\, H \sqcap g_{\underline{\hat{j}}} 
      }}\right), 
\end{eqnarray*}
where $C_{\underline{\hat{i}},\, \underline{\hat{j}}}$ is the number of
elements in the set
\[
\widetilde{C}_{\underline{\hat{i}},\, \underline{\hat{j}}} : = \{(h,
\tilde{h}) \in H \times H : \beta^r_h (\underline{j}) =
\beta^r_{\tilde{ h}} (\underline{\tilde{i}}),\ \beta^r_h
(\underline{i}) = \underline{\hat{i}}\ \mbox{and}\ \beta^r_{\tilde{
    h}} (\underline{\tilde{j}}) = \underline{\hat{j}}\}.
\]
We have used the fact that since the above product is in $N' \cap
M_{2r}$, the coefficient of $\hat{h} [\underline{\hat{i}},\,
\underline{\hat{j}}]^{ev}$ in it is the same as that of
$[\underline{\hat{i}},\, \underline{\hat{j}}]^{ev}$ for all $\hat{h}
\in H$.

On the other hand, by (\ref{Mk}) for the bipartite graph $\star_n$, we
have
\begin{eqnarray*}
  \lefteqn{ \varphi_{2r+1} \left(\sum_{h \in\, H}h  [\underline{i},\,
      \underline{j}]^{ev}\right) \
    \varphi_{2r+1} \left(\sum_{\tilde{h} \in\, H}\tilde{h}
      [\underline{\tilde{i}},\,
      \underline{\tilde{j}}]^{ev}\right) }\\ 
  & = & \sum_{g\, \in\, G} g \left( \vcenter{\small \xymatrix @-2.5pc{ 
        & H g_{{i}_r}, &\, \cdots, &\, H \sqcap g_{\underline{i}} \\
        H & & & \\
        & H g_{{j}_r}, &\, \cdots, &\, H \sqcap g_{\underline{j}}
      }} \right)  
  \sum_{\tilde{g}\, \in\, G} \tilde{g} \left( 
    \vcenter{\small \xymatrix @-2.5pc{
        & H g_{\tilde{i}_r}, &\, \cdots, &\, 
        H \sqcap g_{\underline{\tilde{i}}} \\
        H & & & \\
        & H g_{\tilde{j}_r}, &\, \cdots, &\, H \sqcap g_{\underline{\tilde{j}}} 
      }} \right) 
\end{eqnarray*}

\begin{eqnarray*}
  & = & \sum_{g,\, \tilde{g}\, \in G} \delta^{g\, (H,\, H g_{j_r},\, \cdots,\, 
    H \sqcap g_{\underline{j}})}_{ 
    \tilde{g}\, (H,\, H g_{\tilde{i}_r},\, \cdots,\, H 
    \sqcap g_{\underline{\tilde{i}}})} 
  \left( \vcenter{\small \xymatrix @-2.5pc{
        & g\cdot  H g_{{i}_r}, &\, \cdots, &\ g\cdot H 
        \sqcap g_{\underline{{i}}} \\
        g\cdot H & & & \\
        &\tilde{g}\cdot H g_{\tilde{j}_r}, &\, \cdots, &\ \tilde{g}\cdot H 
        \sqcap g_{\underline{\tilde{j}}} 
      }} \right)\\
  & = & \sum_{\hat{g}\, \in\, G,\ ( \underline{\hat{i}},\, \underline{\hat{j}})\, 
    \in\, H \backslash ( I^r \times I^r)} 
  D_{\underline{\hat{i}},\, \underline{\hat{j}}}\
  \hat{g}  \left(\vcenter{\small \xymatrix @-2.5pc{
        &  H g_{\hat{i}_r}, &\, \cdots, &\,  H \sqcap g_{\underline{\hat{i}}} \\
        H & & & \\
        & H g_{\hat{j}_r}, &\, \cdots, &\,  H \sqcap g_{\underline{\hat{j}}} 
      }} \right),
\end{eqnarray*}
where $D_{\underline{\hat{i}},\, \underline{\hat{j}}}$ is the number
of elements in the set
\[
\widetilde{D}_{\underline{\hat{i}},\, \underline{\hat{j}}} :=\left\{
  \begin{array}{l}
    (g, \tilde{g}) \in G \times G : \\ 
    \!   \left\{ \begin{array}{l} g (H, H g_{j_r},
        \cdots,\, H \sqcap g_{\underline{j}}) = \tilde{g} ( H, H
        g_{\tilde{i}_r}, \cdots,\, H \sqcap g_{\underline{\tilde{i}}}),\\
        g (H, H g_{i_r}, \cdots,\, H \sqcap g_{\underline{i}}) = (
        H, H g_{\hat{i}_r}, \cdots,\, H \sqcap g_{\underline{\hat{i}}})\
        \mbox{and}\\
        \tilde{g} (H, H g_{\tilde{j}_r}, \cdots,\, H \sqcap
        g_{\underline{\tilde{j}}}) = ( H, H g_{\hat{j}_r}, \cdots,\, H \sqcap
        g_{\underline{\hat{j}}})
\end{array} \right.
\end{array}
\! \right\},
\]
and the coefficients are constant on each orbit as in the former case.

It can be seen that, for each $ ( \underline{\hat{i}},\,
\underline{\hat{j}}) \in H \backslash ( I^r \times I^r)$, the sets
$\widetilde{C}_{\underline{\hat{i}},\, \underline{\hat{j}}}$ and
$\widetilde{D}_{\underline{\hat{i}},\, \underline{\hat{j}}}$ are same.

Thus we conclude that $ \varphi_{2r+1} \circ Z^P_{M_{2r+1}} = Z^{Q
}_{M_{2r+1}}\circ (\varphi_{2r+1} \otimes \varphi_{2r+1})$.
\end{pf}
\newline

Note that $Q_{1, k} = [P(\star_n)^G]_{1, k} = P_{1, k}(\star_n)^G,\
\forall\, k \geq 1$.

\begin{lem}\label{conditional-lemma} With running notations, we have
\begin{eqnarray}
  \varphi_{k} \circ Z^P_{E_{k}^{k+1}} & = & Z^Q_{E_{k+1}^k}\circ \varphi_{k+1},\ 
  \mbox{and} \label{con-1} \\
  \varphi_{k+1} \circ Z^P_{(E')_{k+1}^{k+1}} & = &  Z^{Q}_{(E')_{k+1}^{k+1}}\circ 
  \varphi_{k+1},\ \forall\, k \in Col.
  \label{con-2}
\end{eqnarray}
\end{lem}
\begin{pf}
  First note that the maps $\varphi_k, k \in Col$ are all
  $*$-preserving, where the $*$-structure on $Q$ is given as in the
  last paragraph of $\S 4.1$; and, by Lemma
  \ref{multiplication-lemma}, they are algebra homomorphisms as well.

  We next show that these maps preserve the traces as well, where $Q $
  is equipped with the global pictorial trace as given in Lemma
  \ref{p-n-star-spherical}. Then by the fact that they are
  $*$-preserving algebra isomorphisms, it follows that the global
  pictorial trace on $Q$ is in fact faithful. Thus (\ref{con-1}) holds
  by Lemma \ref{conditional-expectation-fact}, and also by the same
  result, (\ref{con-2}) will hold once we show that $\varphi_{k}
  (P_{1,\, k}) = Q_{1,\, k} := Image\, (Z^{Q}_{(E')_{k}^{k}}),
  \forall\, k \geq 1$.
  
  We calculate the trace on $Q_k$ for odd $k$, say $k =
  2r +1$ for some $r \geq 1$. 

  \noindent In the planar algebra $P(\star_n)$, given
  $\underline{i},\, \underline{j} \in I^r$ and $x, y \in I$, we have

  \[
  \begin{array}{rl}
    Z^{P(\star_n)}_{E^{0_+}_{2r+1}}
    \left( \left(\vcenter{\small \xymatrix @-2.5pc{
            & H g_{i_1}, &\, \cdots, &\, H  g_{i_r} \\
            H g_x & & & \\
            & H g_{j_1}, &\, \cdots, &\, H g_{j_r} }}\right)
    \right)_{_{H g_y}} & = 
    \! \vcenter{ \mbox{
        \psfrag{x}{$^{_{H g_x}}$}
        \psfrag{s}{$^{_{H g_y}}$}
        \psfrag{i1}{$^{_{H g_{i_1}}}$}
        \psfrag{ir}{$^{_{H g_{i_r}}}$}
        \psfrag{j1}{$^{_{H g_{j_1}}}$}
        \psfrag{jr}{$^{_{H g_{j_r}}}$}
        \psfrag{*}{$^{_{\ast}}$}
      \includegraphics[scale=0.55]{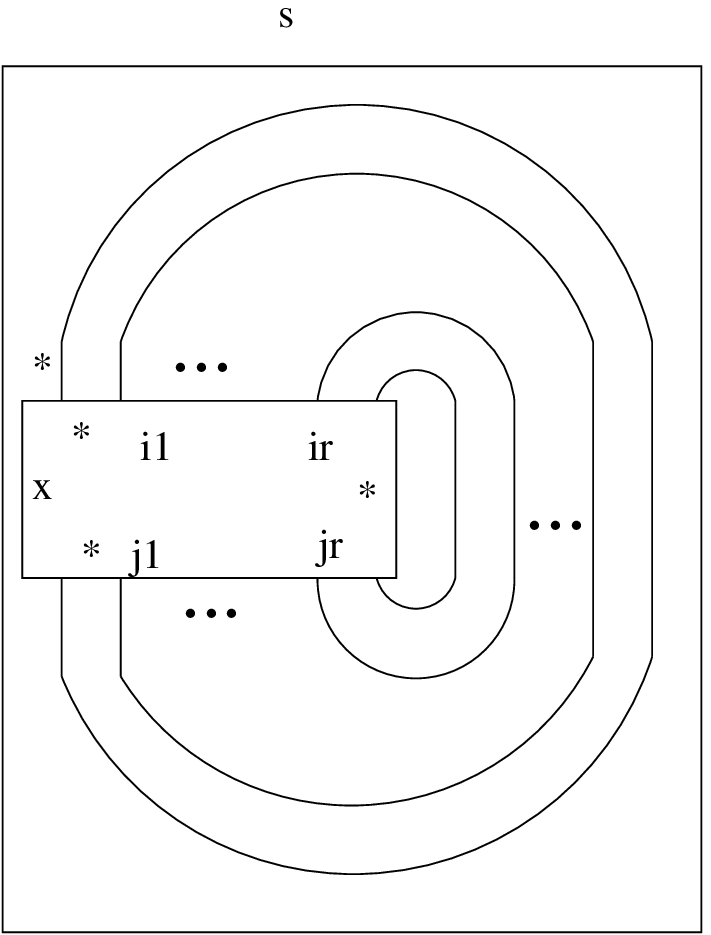}  }} \\
  & =\quad  \delta^{\underline{i}}_{\underline{j}}\ \delta^{x}_{y}\ 
  \mu_{\ast}^2. 
\end{array}
\] 
Thus the pictorial trace on $P_{2r+1}(\star_n)^G$ is given by
    \begin{eqnarray*}
    tr_{2r+1} \left(\sum_{g \in G} g\left(\vcenter{\small \xymatrix
          @-2.5pc{
            & H g_{i_1}, &\, \cdots, &\, H  g_{i_r} \\
            H g_x & & & \\
            & H g_{j_1}, &\, \cdots, &\, H g_{j_r} }}\right) \right) & = &
    \frac{1}{(\sqrt{n})^{2r+1}}\, \mu_{\ast}^2 \,
    \delta^{\underline{i}}_{\underline{j}}\ |\mbox{Iso}_G (H g_x)| \\
& =  &
    \frac{|H|}{n^{r}}\, \delta^{\underline{i}}_{\underline{j}},
  \end{eqnarray*}
for all $x \in I$ and $ \underline{i},\, \underline{j} \in I^r $,
where we have used the fact that, $G$ acting transitively on $H
\backslash G$, $\mbox{Iso}_G (H g_x) \cong \mbox{Iso}_G (H ) = H$, for
all $x \in I$.

On the other hand, for each $(\underline{i},\, \underline{j}) \in H
\backslash (I^r \times I^r) $, we have

$tr_{M_{I^r}(M)} \left([\underline{i},\, \underline{j}]^{ev} \right) =
\delta^{\underline{i}}_{\underline{j}}\, \frac{1}{n^r} $; so that
 \[
 tr_{M_{I^r}(M)} \left(\sum_{h \in H}h [\underline{i},\,
 \underline{j}]^{ev} \right) = \sum_{h \in H}\,
 \delta^{\beta^r_h(\underline{i}) }_{\beta^r_h(\underline{j})}\,
 \frac{1}{n^r} = \frac{|H|}{n^r}\, \delta^{\underline{i}}_{\underline{j}},
\]
and thus
\begin{eqnarray*}
  tr_{2r+1}\left(\varphi_{2r+1} \left(\sum_{h \in H} h
      [\underline{i},\, \underline{j}]^{ev} \right)\right) &=& tr_{2r+1}
  \left( \sum_{g \in G} g\left(\vcenter{\small \xymatrix @-2.5pc{
          & H g_{i_1}, &\, \cdots, &\, H  g_{i_r} \\
          H  & & & \\
          & H g_{j_1}, &\, \cdots, &\, H g_{j_r} }}\right)\right)\\
  & = &
  \frac{|H|}{n^r}\, \delta^{\underline{i}}_{\underline{j}}.
\end{eqnarray*}
This shows that $\varphi_{2r+1}$ preserves trace. With exactly similar
calculations we can show that $\varphi_{k}$ preserves trace for even
$k$ as well.

It only remains to show that $\varphi_{k} (P_{1, k}) = Q_{1, k}$ for
all $k \geq 1$. Again, there is nothing to prove for $k = 1$. We prove
the assertion for odd $k$, say $k = 2r+1$, for
some $r \geq 1$. Note that \\
$\left\{\sum_{g \in G} \left( \vcenter{\small \xymatrix @-2.5pc{ & H
        g_{i_r}, &\, H g_{i_{r-1}} g_{i_r}, &\, \cdots, &\,
        H \sqcap g_{\underline{i}} \\
        H  & & & \\
        & H g_{j_r},&\, H g_{j_{r-1}} g_{j_r},&\, \cdots, &\, H\sqcap
        g_{\underline{j}} }}\right) : (\underline{i},\, \underline{j}
  ) \in H \backslash (I^r \times I^r) \right\}$ forms a basis for
$Q_{2r+1}$.

By (\ref{E'k}) for the bipartite graph $\star_n$, for each
$(\underline{i},\, \underline{j} ) \in H \backslash (I^r \times I^r),$
we have
\begin{eqnarray*}
  \lefteqn{  Z^{Q}_{(E')_{k}^{k} } \left( \sum_{g \in G} g  \left(
        \vcenter{\small \xymatrix @-2.5pc{
            & H g_{i_r},&\,  Hg_{i_{r-1}} g_{i_r}, &\, \cdots, &\,
            H \sqcap g_{\underline{i}} \\
            H  & & & \\
            & H g_{j_r}, &\,  Hg_{j_{r-1}} g_{j_r},  &\, \cdots, &\,
            H \sqcap g_{\underline{j}} }}
      \right)\right)} \\
  & = & \frac{1}{\sqrt{n}} \sum_{z \in I,\, g \in G} \left(
    \vcenter{\small \xymatrix @-2.5pc{
        & g \cdot H g_{i_r}, &\, g \cdot Hg_{i_{r-1}} g_{i_r}, &\,
        \cdots, &\, g \cdot H \sqcap  g_{\underline{i}} \\
        H g_z & & & \\
        & g \cdot H g_{j_r}, &\, g \cdot Hg_{j_{r-1}} g_{j_r}, &\, \cdots, &\,
        g \cdot H \sqcap g_{\underline{j}} }}\right)\\
  & \qquad = &\frac{1}{\sqrt{n}} \sum_{x \in I,\, g \in G} g \left( 
    \vcenter{\small  \xymatrix @-2.5pc{
        & H g_{i_r},&\,  Hg_{i_{r-1}} g_{i_r},  &\, \cdots, &\, 
        H \sqcap g_{\underline{i}} \\
        H g_x & & & \\
        & H g_{j_r},&\,  Hg_{j_{r-1}} g_{j_r},  &\, \cdots, &\, 
        H \sqcap g_{\underline{j}} }}\right);
\end{eqnarray*}
and these elements generate $Q_{1,\, 2r+1} $ as a vector space.

Further, for each $ (\underline{i},\, \underline{j} ) \in H \backslash
(I^r \times I^r) $, $\sum_{g \in G} g [\underline{i},\,
\underline{j}]^{ev} \in M' \cap M_{2r}$ and
\begin{eqnarray*}
  \sum_{g \in G} g [\underline{i},\, \underline{j}]^{ev} & = &
  \sum_{x \in I,\, h \in H} 
  h g_x [\underline{i},\, \underline{j}]^{ev}\\
  \stackrel{\varphi_{2r+1}}{\longmapsto} & & \sum_{x \in I,\, g \in G} g
  \left(
    \vcenter{\small \xymatrix @-2.5pc{
        & H g_{i_r}g_x^{-1}, &\,  Hg_{i_{r-1}} g_{i_r}g_x^{-1}, &\, \cdots,
        &\, H (\sqcap g_{\underline{i}}) g_x^{-1} \\
        H  & & & \\
        & H g_{j_r}g_x^{-1}, &\,  Hg_{j_{r-1}} g_{j_r}g_x^{-1}, &\, \cdots, &\, 
        H (\sqcap g_{\underline{j}}) g_x^{-1} }}\right) \\
  & = & \sum_{x \in I,\, g \in G} g g_x
  \left(
    \vcenter{\small \xymatrix @-2.5pc{
        & H g_{i_r},&\,  Hg_{i_{r-1}} g_{i_r}, &\, \cdots, &\, H  
        \sqcap g_{\underline{i}} \\
        H g_x & & & \\
        & H g_{j_r},&\,  Hg_{j_{r-1}} g_{j_r}, &\, \cdots, &\, H
        \sqcap g_{\underline{j}} }}\right)\\
  & = & \sum_{x \in I,\, \tilde{g} \in G} \tilde{g} 
  \left(
    \vcenter{\small \xymatrix @-2.5pc{
        & H g_{i_r},&\,  Hg_{i_{r-1}} g_{i_r}, &\, \cdots, &\, H
        \sqcap  g_{\underline{i}} \\
        H g_x & & & \\
        & H g_{j_r}, &\,  Hg_{j_{r-1}} g_{j_r},&\, \cdots, &\, 
        H \sqcap g_{\underline{j}} }}\right)
\end{eqnarray*}
Thus the elements of above generating set for $Q_{1, 2r+1}$ are in the
space $\varphi_{2r+1} (M' \cap M_{2r})$, and we know that $P_{1, 2r+1}
= E_{M' \cap M_{2r}}(N' \cap M_{2r}) = M' \cap M_{2r}$.  This proves
our second claim.
\end{pf}
\newline

\begin{lem}\label{jones-projection-lemma}
  $ ( \varphi_{k+1} \circ Z^P_{\mathcal{E}^{k+1}} ) =
  Z^{Q}_{\mathcal{E}^{k+1}},\ \forall \, k \geq 1.  $
\end{lem}
\begin{pf}
  We first prove the assertion for even $k$.

  We need to consider some $H$ and $G$ invariant subsets of $I^k
  \times I^k$ and $(H \backslash G)^{2k+1}$ for all $k \geq 1$. For
  each $k \geq 1$, we set $W_k = \{ (\underline{i},\, \underline{j})
  \in I^k \times I^k : \underline{i}_{[2} = \underline{j}_{[2}\} $,
\[
\begin{array}{l}
  F_{2k+1}\!  =  \left\{\! \left( \vcenter{\small \xymatrix @-2.5pc{
          & H g_{i_1}, &\, \cdots, &\, H  g_{i_k} \\
          H g_x & & & \\
          & H g_{j_1}, &\, \cdots, &\, H g_{j_k} }}\right) : x \in I,\,
    \underline{i}, \underline{j} \in I^{k}\ \mbox{with}\ 
    \underline{i}_{\, k-1]} =    \underline{j}_{\, k-1]}\! \right\}\\
  \mbox{and}   \\
  G_{2k+1}\!  =  \left\{\! \left( \vcenter{\small \xymatrix @-2.5pc{
          & H g_{i_1}, &\, \cdots, &\, H  g_{i_k} \\
          H  & & & \\
          & H g_{j_1}, &\, \cdots, &\, H g_{j_k} }}\right) :
    \underline{i}, \underline{j} \in I^{k}\ \mbox{with}\
    \underline{i}_{\, k-1]} = \underline{j}_{\, k-1]}\! \right\}.
\end{array}
\]
Then $W_k$ (resp., $F_{2k+1}$) is $G$ invariant under the diagonal
$\beta^k$(resp., $\beta^1$)-action, and $G_{2k+1} \subset F_{2k+1} $
is $H$ invariant under the restricted action.

Further, the correspondence
\[
W_k \ni (\underline{i},\, \underline{j}) \mapsto \left(
  \vcenter{\small \xymatrix @-2.5pc{ & H g_{i_k},&\, Hg_{i_{k-1}}
      g_{i_k}, &\, \cdots, &\, H   \sqcap g_{\underline{i}} \\
      H  & & & \\
      & H g_{j_k},&\, Hg_{j_{k-1}} g_{j_k}, &\, \cdots, &\, H \sqcap
      g_{\underline{j}} }}\right) \in G_{2k+1}
\]
is an $H$-bijection. So $\left\{\! \left( \vcenter{\small \xymatrix
      @-2.5pc{ & H g_{i_k},&\, Hg_{i_{k-1}}
        g_{i_k}, &\, \cdots, &\, H \sqcap g_{\underline{i}} \\
        H  & & & \\
        & H g_{j_k}, &\, Hg_{j_{k-1}} g_{j_k}, &\, \cdots, &\, H
        \sqcap g_{\underline{j}} }} \right) : (\underline{i},\,
  \underline{j}) \in H \backslash W_k \right\}$ is a set of
representatives of $H$ orbits of $G_{2k+1}$.  Clearly this is also a
set of representatives of $G$ orbits of $F_{2k+1}$.

Now we note that, by (\ref{Ek+1}) for the graph $\star_n$, we have
\[
Z^{Q}_{\mathcal{E}^{2k+1}} (1) = \frac{1}{\sqrt{n}}
\sum_{x, y, z \in I,\, \underline{i} \in I^{k-1} }
\left( \vcenter{\small \xymatrix @-2.5pc{
        & H g_{i_1}, &\, \cdots, &\, H  g_{i_{k-1}}, &\, H  g_y \\
        H g_x  & & & \\
        & H g_{i_1}, &\, \cdots, &\, H g_{i_{k-1}}, &\, H  g_z }}\right).
\]
Recall, from Theorem \ref{jones}, that $Z^{P}_{\mathcal{E}^{2k+1}} (1)
= \sqrt{n}\, \widetilde{e}_{2k}$, where the Jones projection
$\widetilde{e}_{2k} \in M_{I^k}(M)$ is given, as in Corollary
\ref{model-facts}, by
\[
(\widetilde{e}_{2k})_{\underline{i},\, \underline{j}} = n^{-1}\,
\delta^{\underline{i}_{[2}}_{\underline{j}_{[2}}\, g_{i_1}
g_{j_1}^{-1},\ \forall\ \underline{i}, \underline{j} \in I^k.
\]
Thus 
\begin{eqnarray*}
  \widetilde{e}_{2k} & = &  \frac{1}{n} \sum_{(\underline{i},\,
    \underline{j})\, \in\, W_k} [\underline{i},\, \underline{j}]^{ev}\\
  & = &
  \frac{1}{n} \sum_{h\, \in\, H, \atop{ (\underline{i},\, \underline{j})\,
      \in\, H \backslash W_k}}
  \frac{1}{|\mbox{Iso}_H (\underline{i},\, \underline{j})|}\ h
  [\underline{i},\, \underline{j}]^{ev} \\
  \stackrel{\varphi_{2k+1}}{\longmapsto} & & \frac{1}{n} \sum_{g\, \in\, G,
    \atop{ (\underline{i},\, \underline{j})
      \in\, H \backslash W_k}}
  \frac{1}{|\mbox{Iso}_H (\underline{i},\, \underline{j})|}\ g
  \left(
    \vcenter{\small \xymatrix @-2.5pc{
        & H g_{i_k}, &\, \cdots, &\, H \sqcap g_{\underline{i}} \\
        H  & & & \\
        & H g_{j_k}, &\, \cdots, &\, H \sqcap g_{\underline{j}}
      }}\right)\\
  & = & \frac{1}{n} \sum_{g\, \in\, G,
    \atop{ (\underline{i},\, \underline{j})\, \in\, H \backslash W_k}}
  \frac{1}{\left|\mbox{Iso}_G \left(\left(
        \vcenter{\small \xymatrix @-2.5pc{
            & H g_{i_k}, &\, \cdots, &\, H \sqcap g_{\underline{i}} \\
            H  & & & \\
            & H g_{j_k}, &\, \cdots, &\, H \sqcap g_{\underline{j}}
          }}\right) \right) \right|}\ \times \\
  & & \hspace*{50mm} \ g
  \left(
    \vcenter{\small \xymatrix @-2.5pc{
        & H g_{i_k}, &\, \cdots, &\, H \sqcap g_{\underline{i}} \\
        H  & & & \\
        & H g_{j_k}, &\, \cdots, &\, H \sqcap g_{\underline{j}}
      }}\right)\\
  & = & \frac{1}{n}
  \sum_{x,\, y,\, z\, \in\, I,\ \underline{i}\, \in\, I^{k-1} }
  \left( \vcenter{\small \xymatrix @-2.5pc{
        & H g_{i_1}, &\, \cdots, &\, H  g_{i_{k-1}}, &\, H  g_y \\
        H g_x  & & & \\
        & H g_{i_1}, &\, \cdots, &\, H g_{i_{k-1}}, &\, H  g_z }}\right),
\end{eqnarray*}
where we have used the fact that 
\[
\mbox{Iso}_H ((\underline{i},\, \underline{j})) = \mbox{Iso}_G
\left(\left( \vcenter{\small \xymatrix @-2.5pc{
        & H g_{i_k}, &\, \cdots, &\, H \sqcap g_{\underline{i}} \\
        H  & & & \\
        & H g_{j_k}, &\, \cdots, &\, H \sqcap g_{\underline{j}}
      }}\right) \right).
\]
This shows that $ Z^{Q}_{\mathcal{E}^{2k+1}} =
\varphi_{2k+1} \circ Z^{P}_{\mathcal{E}^{2k+1}} $.

The case for odd $k$ follows on similar lines by taking analogues of
$W_k$, $F_{2k+1}$ and $G_{2k+1}$ to be the sets $Z_k$, $F_{2k}$,
$G_{2k}$, respectively, which are given by
\[
  Z_k  =  \{(\underline{i},\, \underline{j} ) \in I^k \times I^k : 
  \underline{i} = \underline{j}, i_1 = 1\},
\]
\[
\begin{array}{l}
  F_{2k}\!  = \!  \left\{\!\!
    \left(\! \! \vcenter{\small \xymatrix @-2.5pc{
          & H g_{i_1}, &\, \cdots, &\, H  g_{i_{k-1}}& \\
          H g_x & & & & H g_y \\
          & H g_{i_1}, &\, \cdots, &\, H g_{i_{k-1}}& }}\!\! \right)\!\!:
    x, y \in I,\,
    \underline{i} \in I^{k-1}\ \mbox{with}\ y = i_{k-1}\!\!
  \right\}\\
  \mbox{and} \\
  G_{2k} \! =\!  \left\{ \!
    \left(\! \vcenter{\small \xymatrix @-2.5pc{
          & H g_{i_1}, &\, \cdots, &\, H  g_{i_{k-1}}& \\
          H  & & & & H g_y \\
          & H g_{i_1}, &\, \cdots, &\, H g_{i_{k-1}}& }}\! \right)\!\!: 
    y \in I,\, \underline{i} \in I^{k-1}\ \mbox{with}\ y = i_{k-1} \!
  \right\}.
\end{array}
\]
\end{pf}
\newline 

Finally, all the calculative work being done, we collect some of the
lemmas proved above to give a complete proof of the main theorem.
\begin{thm}\label{main-theorem}
  Given a finite group $G$, a subgroup $H$ of index, say $n$, and an
  outer action $\alpha$ of $G$ on the hyperfinite $II_1$-factor $R$,
  the planar algebra of the subgroup-subfactor ${R \rtimes H \subset R
    \rtimes G}$ is isomorphic to the $G$-invariant planar subalgebra
  of $P (\star_n)$, i.e.,
\[
P^{R \rtimes H \subset R \rtimes G} \, \cong \, P(\star_n)^G.
\]
\end{thm}
\begin{pf} 
  Let $\varphi_k,\, k \in Col$ be the maps defined as in equations
  (\ref{varphi-0-pm}$-$\ref{varphi-k-even}). We claim that
\[
\{\varphi_k : k \in Col\} =: \varphi : P^{R \rtimes H \subset R
  \rtimes G} \rightarrow P(\star_n)^G
\]
is a planar algebra isomorphism.

We already know, by Lemma \ref{varphi-isomorphism}, that the maps
$\varphi_k, k \in Col$ are all linear isomorphisms. Thus what remains
to be shown is that $\varphi$ is a planar algebra morphism.

Let $\mathcal{T}$ be the collection of coloured tangles $T$ which
commute with $\varphi$.  Then, by Theorem \ref{gen-tangles}, it is
enough to show that $\mathcal{T}$ is closed under composition of
tangles, whenever it makes sense, and that it contains the generating
set of tangles $\mathcal{G}_0 = \{1^{0_{\pm}}\} \cup \{ E^k_{k+1},
M_k, I^{k+1}_{k}:\, k \in\, Col \} \cup \{\mathcal{E}^{k+1}, (E')^k_k
: k \geq 1\}.$

The first assertion follows from Lemma \ref{closedness-of-tangles}.
Then, it readily follows from definitions that $\varphi_{0_{\pm}}
\circ Z^P_{1^{0_{\pm}}} = Z^{P(\star_n)^G}_{1^{0_{\pm}}}$. Thus
$\{1^{0_+}, 1^{0_-}\} \subset \mathcal{T}$.  Finally, Lemmas
\ref{inclusion-lemma}$-$\ref{jones-projection-lemma} show that the
collection $ \{I^{k+1}_{k}, M_k :\, k \in\, Col \} \cup \{E^k_{k+1},
(E')^{k+1}_{k+1} : k \in Col \}\cup \{\mathcal{E}^{k+1} : k \geq 1\} $
is also contained in $\mathcal{T}$.  Thus $\mathcal{G}_0 \subset
\mathcal{T}$.

This completes the proof.
\end{pf}
\newline

What follows immediately is the following fact, which, however,
already exists in literature.
\begin{cor}
  Given a finite group $G$ and a subgroup $H$, the subgroup-subfactor
  $R \rtimes_{\alpha \slash H} H \subset R \rtimes_{\alpha} G$ does
  not depend upon the outer action $\alpha$ of $G$ on the hyperfinite
  $II_1$-factor $R$.
\end{cor}
\begin{pf}
  It follows from Theorem \ref{main-theorem} that the planar algebra
  of the subgroup-subfactor $R \rtimes_{\alpha \slash H} H \subset R
  \rtimes_{\alpha} G$ is independent of the outer action $\alpha$.  In
  particular, it follows that the standard invariant of the subfactor
  is independent of $\alpha$.  $R \rtimes_{\alpha} G$ being
  hyperfinite, the standard invariant of $R \rtimes_{\alpha \slash H}
  H \subset R \rtimes_{\alpha} G$ is a complete invariant ({\em c.f.}
  \cite{Pop94}). Thus the subgroup-subfactor $R \rtimes_{\alpha \slash
    H} H \subset R \rtimes_{\alpha} G$ is independent of the outer
  action $\alpha$.
\end{pf}
\newline

By repeated applications of the `Not Burnside's Lemma', the
coefficients in the Poincar$\acute{\mbox{e}}$ series of the
subgroup-subfactor are given by:

\begin{cor}\label{dimensions}For each $ k \geq 1,$
  \begin{equation}
dim\ P^{ R \rtimes H \subset R \rtimes G}_k = \frac{1}{|G|} \sum_{C\,
    \in\,  \mathcal{C}_{G}} |C| \left( \frac{|C \cap H|
      |G|}{|C||H|}\right)^k,
\end{equation}
where $\mathcal{C}_{G}$ is the set of conjugacy classes of $G$.
\end{cor}

As promised above, we note that since the automorphism group of the
bipartite graph $\star_n$ is $S_n$, we have the following:
\begin{cor}\label{sandwich}
  Given any pair of finite groups $H \subset G$ with index $n$, the
  planar algebra $P^{R \rtimes H \subset R \rtimes G}$ is a planar
  subalgebra of the planar algebra of the bipartite graph $\star_n$
  and contains the planar algebra $P(\star_n)^{S_n} \cong P^{R \rtimes
    S_{n-1} \subset R \rtimes S_n}$.
\end{cor}

We can also identify the planar algebra of the fixed subalgebra with
the $G$-invariant planar subalgebra of the planar algebra of the flip
of $\star_n$.
\begin{cor}\label{dual-theorem}
  Given a finite group $G$, a subgroup $H$ of index, say $n$, and an
  outer action $\alpha$ of $G$ on the hyperfinite $II_1$-factor $R$,
  the planar algebra $P^{R^G \subset R^H}$ is isomorphic to the
  $G$-invariant planar subalgebra of $P (\overline{\star_n})$, i.e.,
\[
P^{R^G \subset R^H} \, \cong \, P (\overline{\star_n})^G,
\]
where $\overline{\star_n}$ is the flip of the bipartite graph
$\star_n$. In particular, the subfactor $R^G \subset R^H$ is
independent of the outer action $\alpha$.
\end{cor}
\begin{pf}
  Recall, from Proposition \ref{fixed-algebra-subfactor}, that if $N
  := R \rtimes_{\alpha \slash H} H$ and $M := R \rtimes_{\alpha} G$,
  then $P^{R^G \subset R^H} \cong P^{M \subset M_1}$, where $M_1$ is
  the $II_1$-factor obtained by the basic construction of the
  subfactor $N \subset M$. Further, by \cite[Proposition 4.17]{KS04},
  $P^{M \subset M_1}$ is isomorphic to the dual planar algebra $^-
  P^{N \subset M}$. Thus, by Theorems \ref{gp-action-p-gamma} and
  \ref{main-theorem}, we conclude that
  \[
  P^{R^G \subset R^H} \, \cong \, P (\overline{\star_n})^G.
  \]
This shows that the planar algebra $P^{R^G \subset R^H}$ is
independent of the outer action $\alpha$.  Further, $R^H$, being a
$II_1$-factor sitting in $R$, is itself hyperfinite - see
\cite{Con76}.  Thus, as in the preceding corollary, the subfactor $R^G
\subset R^H$ is independent of the action $\alpha$.
\end{pf}

\section*{Acknowledgements}
The author would like to thank his advisor Prof.~V.~S.~Sunder and
Prof.~Vijay Kodiyalam for their invaluable guidance and support
throughout his stay at IMSc as a graduate student. This paper is a
result of many long discussions that the author had with them, and it
would not have come into existence without their input and guidance.

\nocite{Jon83}
\bibliographystyle{plain}
\bibliography{pln-alg-sbgp-sbfctr}

\end{document}